\documentclass[natbib,11pt,a4]{article}
\usepackage{subfigure}% Support for small, `sub' figures and tables
\RequirePackage{fix-cm}

\usepackage{amsmath}
\usepackage{amssymb}
\usepackage{theorem}
\usepackage{array}
 \usepackage{subfigure}
   \usepackage{xspace}
\usepackage{epsfig}
\usepackage{hhline}
\usepackage{bbm}
\usepackage{natbib}
\usepackage{nonfloat}
\usepackage[utf8]{inputenc}    % pour les accents 

\theoremstyle{plain}
\newtheorem{theorem}{Theorem}[section]

\newtheorem{lemma}[theorem]{Lemma}
\newtheorem{proposition}[theorem]{Proposition}

\theoremstyle{remark}
\newtheorem{remark}[theorem]{Remark}

\theoremstyle{definition}

\def\pth#1{\left(#1\right)}
\def\acc#1{\left\{#1\right\}}
\def\cro#1{\left[#1\right]}

\def\efo{\textrm{\mathversion{bold}$\mathbf{\phi^0}$\mathversion{normal}}}
\def\uuu{\textrm{\mathversion{bold}$\mathbf{1}$\mathversion{normal}}}
\def\oo{\textrm{\mathversion{bold}$\mathbf{0}$\mathversion{normal}}}
  
\def\el{\textrm{\mathversion{bold}$\mathbf{\lambda}$\mathversion{normal}}}  
 
\def\eU{\textrm{\mathversion{bold}$\mathbf{\Upsilon}$\mathversion{normal}}} 
\def\eo{\textrm{\mathversion{bold}$\mathbf{\omega}$\mathversion{normal}}} 

\def\eE{I\!\!E}
\def\eP{I\!\!P}
\def\e1{1\!\!1}

\def\eeX{\mathbb{X}}
\def\hh{ \hspace*{0.5cm}}

\def\ef{\textrm{\mathversion{bold}$\mathbf{\phi}$\mathversion{normal}}}  
\def\ee1{\textrm{\mathversion{bold}$\mathbf{\varepsilon}$\mathversion{normal}}}

\def\XX{\textrm{\mathversion{bold}$\mathbf{X}$\mathversion{normal}}}
\def\xx{\textrm{\mathversion{bold}$\mathbf{x}$\mathversion{normal}}}
\def\uu{\textrm{\mathversion{bold}$\mathbf{u}$\mathversion{normal}}}

\newcommand{\R}{\mathbb{R}}
\newcommand{\N}{\mathbb{N}}

\def\argmin{\mathop{\mathrm{arg\,min}}} 

\begin{document}
%\bibliographystyle{plan} 
%--------------------------------------------------

\title {Adaptive LASSO model selection in a multiphase quantile regression }
%\author{Gabriela Ciuperca         }


%\institute{ G.Ciuperca  \at
%Universit\'e de Lyon, Universit\'e Lyon 1, CNRS, UMR 5208, Institut Camille Jordan, 
 %            Bat.  Braconnier, 43, blvd du 11 novembre 1918, F - 69622 Villeurbanne Cedex, France \\
   %                      \email{Gabriela.Ciuperca@univ-lyon1.fr}          
  %      }
\author{Gabriela Ciuperca $^{\ast}$\thanks{$^\ast$Corresponding author. Email: Gabriela.Ciuperca@univ-lyon1.fr
\vspace{6pt}} \\\vspace{6pt}  {\em{Universit\'e de Lyon, Universit\'e Lyon 1, CNRS, UMR 5208, Institut Camille Jordan,}}\\
{\em{Bat.  Braconnier, 43, blvd du 11 novembre 1918, F - 69622 Villeurbanne Cedex, France}}\\ }

\maketitle

%%%%%%%%%%%%%%%%%%%%%%%%%%%%%%%%%%%%%%%%%%%%%%%%%%%%%%%%%%%%%%%%%%%%%%%%%%%%%%%%%%%%%%%%%%%%%%%%%%%%%%%%%%%%%%%%%%%%%%%%%%%%

\begin{abstract}
We propose a general adaptive LASSO method for a quantile regression model. Our method is very interesting when we know nothing about the first two moments of the model error. We first prove that the obtained estimators satisfy the oracle properties, which involves the relevant variable selection without using hypothesis test. Next, we study the proposed method when the (multiphase) model changes to unknown observations called change-points. Convergence rates of the change-points and of the regression parameters estimators in each phase are found. The sparsity of the adaptive LASSO quantile estimators of the regression parameters is not affected by the change-points estimation. If the phases number is unknown, a consistent criterion is proposed. Numerical studies by Monte Carlo simulations show the performance of the proposed method, compared to other existing methods in the literature, for models with a single phase or for multiphase models.  The adaptive LASSO quantile method performs better than known variable selection methods, as the least squared method with adaptive LASSO penalty, $L_1$-method with LASSO-type penalty and quantile method with SCAD penalty.

\end{abstract}

%%%%%%%%%%%%%%%%%%%%%%%%%%%%%%%%%%%%%%%%%%%%%%%%%%%%%%%%%%%%%%%%%%%%%%%%%%%%%%%%%%%%%%%%%%%%%%%%%%%%%%%%%%%%%%%%%%%%%%%%%%%%
\section{Introduction}
The usual case investigated in literature, for the following regression model,
\begin{equation}
\label{eq1}
Y_i=\XX_i^t\ef+\varepsilon_i, \qquad i=1, \cdots, n,
\end{equation}
 is that the errors are assumed to be  homoscedastic, i.e., the errors are assumed to be independent random variables with mean zero and bounded variance.  In such cases,  the regression parameters $\ef$ are estimated  by the method of the least squares (LS). We will call the regression in this case, the LS model. If the assumptions on the first two moments of the errors are not satisfied,  then the LS method is not appropriate, because it can provide bad estimators (biased, with large variance). In the case of errors with zero mean of sign, i.e. $\eE[sgn(\varepsilon_i)]=0$, the  Least Absolute Deviations (LAD) method could be used. But, often in practice, we can not  known if $\eE[sgn(\varepsilon_i)]=0$, then, a generalization can be used by considering the quantile method. Here, we use the notation $sgn(.)$ for the sign function.  \\   
 \hspace*{0.5cm} For a fixed quantile index(level) $\tau \in (0,1)$, the $\tau$th quantile $b_\tau$ of $\varepsilon$ is:
\begin{equation}
\label{eq2}
\tau= \eP[\varepsilon < b_\tau] =F(b_\tau), \qquad b_\tau \in {\cal B} \subseteq \R, 
\end{equation}
where $F$ is the distribution function of the error $\varepsilon$ and ${\cal B}$ is a real set . We use $\varepsilon$ to denote a generic member of the sample $(\varepsilon_i)_{1 \leq i \leq n}$. In order to consider a general case for model (\ref{eq1}), the $\tau$th quantile $b_\tau$ of $\varepsilon$ is supposed unknown. \\
 \hspace*{0.5cm} For the model (\ref{eq1}), $Y$ is the response variable,  $\ef=(\phi_{,1}, \phi_{,2}, \cdots, \phi_{,p}) \in \Gamma \subseteq \R^p$ the regression parameters and $\XX_i=(X_{i1}, \cdots, X_{ip})$ the regressors. Number $p$ of the regressors can be large. More exactly, $p$ does not depend on the observation number $n$, and $p <n$. \\
 \hspace*{0.5cm} The fact that the $\tau$th quantile of $\varepsilon$ is unknown is also often the case in a change-point model (model with multiple phases), model defined in Section 3 by the relation (\ref{eq11}).  \\
 \hspace*{0.5cm} We will first focus on the study of model (\ref{eq1}),  without assuming classical conditions imposed on errors $\varepsilon_i$.   The model parameters $\ef$ will be estimated initially by a method without penalty, taking the quantile-process as objective function.  Afterwards, in order to select the variables, we  propose to add to the quantile-process an adaptive general penalty of LASSO type. In these cases, we will then refer to the model (\ref{eq1}) as quantile regression. \\ 
 \hspace*{0.5cm} Under the metric $d(x,y)=|\arctan x - \arctan y  |$, for $x, y \in \R$, the set $\bar \R \equiv \R \cup \{ - \infty, \infty\}$ is compact set, then, without loss of generality, we suppose that the sets  ${\cal B}$ and $\Gamma$ are compacts. The $\tau$th conditional quantile function of $Y_i$ given $\XX_i=\xx_i$ is $\xx^t_i \efo+b^0_\tau$, with $b^0_\tau$ the $\tau$th quantile of $\varepsilon$.  Then, for model (\ref{eq1}), the unknown parameters, to estimate, knowing $(Y_i, \XX_i)_{1 \leq i \leq n}$ are the $\tau$th quantile $b_\tau$ and the regression parameters $\ef$. Denote by $b^0_\tau$, $\efo =(\phi^0_{,1}, \cdots, \phi^0_{,p})$ their true values, unknown, assumed to be inner points of the sets  ${\cal B}$, $\Gamma$, respectively.  The parameters $b_\tau, \ef$ can be estimated by quantile method, by minimizing the quantile-process:
\begin{equation}
\label{eq3}
(\hat b_n, \hat \ef_n)\equiv \argmin_{(b,\ef)} \sum^n_{i=1} \rho_\tau(Y_i-b-\XX_i^t \ef),
\end{equation}
with the function $\rho_\tau(.): \R \rightarrow \R$ defined by $\rho_\tau(r)=r[ \tau \e1_{r>0}-(1-\tau) \e1_{r \leq 0} ]$. We call $(\hat b_n, \hat \ef_n)$ the quantile estimators of $(b^0_\tau,\efo)$. The components of $\hat \ef_n$ are $(\hat \phi_{n,1}, \cdots, \hat \phi_{n,p})$. By \cite{Koenker:05}  these estimators are strongly convergent and asymptotically normal:
\[
(\hat b_n, \hat \ef_n) \overset{a.s.} {\underset{n \rightarrow \infty}{\longrightarrow}} (b^0_\tau,\efo), 
\]
\[  
\sqrt{n} (\hat \ef_n -\efo) \overset{\cal L} {\underset{n \rightarrow \infty}{\longrightarrow}} {\cal N} \pth{ \textbf{0}, \frac{ \tau(1-\tau)}{f^2(b^0_\tau)}{\eU}^{-1} }, 
\]
with the matrix  $\eU$ defined later in the assumption (A1) and $f$ the density function of $\varepsilon$. \\
 \hspace*{0.5cm} For $\tau=1/2$ and $b^0_\tau=0$ we obtain the median (LAD, $L_1$) regression,  considered for instance by \cite{Babu:89}. For a complete review on quantile regression and on quantile estimators, the reader can see the book of \cite{Koenker:05}. A great advantage of this method is that, compared to the classical estimation methods (least squares or likelihood methods, for example) that are sensitive to outliers, the quantile method provides more robust estimators. Furthermore, the imposed condition  to the error distribution are  relaxed. However, like all other estimators (for example least squares, maximum likelihood estimators) obtained by optimizing an objective function without penalty, the quantile estimators do not satisfy the oracle properties. Recall that the oracle properties are: the zero components of the true parameters are estimated (shrunk) as 0 with probability tending to 1 (also called sparsity property) and the nonzero parameters have an optimal estimation rate (and is asymptotically normal).  Thus, the solution is to consider a penalized objective function. The parameters for which we want to have the selection consistency property are included in the penalty. This type of penalty was introduced by \cite{Tibshirani:96} for the least squares estimation framework, with a  $L_1$ penalty (obtaining thereby LASSO estimators). Nevertheless the LASSO estimator does not always  satisfy the oracle properties. To remove this inconvenience, an adaptive LASSO  estimator   was proposed by \cite{Zou:06} for the  LS regression case. Various procedures, based on LASSO framework, have been proposed and studied during the last few years, in order  to simultaneously estimate the parameters and to select significant regressors. 
The considered models are either with fixed dimension $p$ for $\ef$ or with $p$ depending on the sample size $n $, with $p>n$. Since the number of papers in these areas is very important in recent years,  we mention only some of them. In a median regression,  \cite{Xu:Ying:10}  consider the LASSO-type penalty and  \cite{Wang:12} proposes a $L_1$ penalty with the possibility that  the regressors number is larger than the observation number. Always for a median model, with LASSO penalty, was considered by \cite{Gao:Huang:10}.  In a general quantile regression, \cite{Wu:Liu:09} propose the SCAD penalty, but which is difficult into practice with regard to numerical algorithms and an adaptive LASSO penalty but under the assumption that the $\tau$th quantile $b^0_\tau$ is known and is equal to 0. Alternatively, \cite{He:Wang:Hong:13} estimate the quantile-adaptive model-free screening frameworks using a B-spline approximation. In the paper of \cite{Zou:Yuan:08}, a composite quantile regression is considered with an adaptive LASSO penalty. The paper of \cite{Guo:Tang:Tian:Zhu:13}  proposes an estimation procedure in a semi-parametric additive partial linear models.  \\
\hh In this paper we first propose, for quantile model (\ref{eq1}), a general adaptive LASSO estimator. For this estimator,  we will study the oracle properties and other behaviour properties for the estimators and for the  adaptive quantile objective function.\\
 Afterwards, we will study  for  a multiphase model (i.e.,   a model that changes the shape to unknown observations), if the break estimation affects the oracle properties of the regression parameter estimators.  The main theoretical difficulty of this type of model is that the estimation of the change-point locations and of the parameters in each phase can not be performed simultaneously, but sequentially, in the sense that we first estimate regression parameters for fixed change-points and then, the change-points. Finally, the estimator of the regression parameters is taken as that corresponding to the optimal change-points. This imply that the theoretical study of obtained estimators is very difficult, even assuming that  number of phases is known. \\
 If this number is unknown, an additional difficulty to the model study is added. In this case, a consistent estimation criterion for this number is proposed in this paper.  Since the LASSO techniques are fairly recent, there are not many papers in the literature that address the breaking problem by this estimation method. In the paper of \cite{Ciuperca:13a}, LS model is estimated by   LASSO-type and by  adaptive LASSO techniques. \\
 \hh In  \cite{Ciuperca:13b}, a quantile model with SCAD estimator and a median model with LASSO-type estimator are studied. Apart from the fact that when the quantile model was considered, it was supposed that the $\tau$th error quantile is known (more precisely, it was taken as 0), the SCAD method has also the disadvantage that it is difficult to put into practice with regard to numerical algorithms. Let us also remind the paper of \cite{Lee:Seo:Shin:12}, where a  LS model with a single  change-point is considered, with a LASSO penalty, when the errors have gaussian distribution.     \\
It is important to emphasize that the numerical studies by Monte Carlo  simulations confirm the superiority of the adaptive LASSO quantile estimators (in terms of bias, precision, of identification of the true zeros, especially in the case when moments of errors don't exist or when its median is different to zero), in comparison  to other variable selection estimators. In a multiphase model, if the changes are due only to error distribution (regression parameters remain the same) the adaptive quantile method also gives the best results.   \\ 
\hh We give some general notations. Throughout the paper, $C$ denotes a positives generic constant not dependent on $n$ which may take different values in different formula or even in different parts of the same formula. For a vector $\mathbf{v}=(v_1, \cdots, v_p)$ let us denote $| \mathbf{v} | =(|v_1|, \cdots, |v_p|)$. On the other hand, $\|\mathbf{v} \|_2$ is the Euclidean norm and $\|\mathbf{v} \|_1= \sum^p_{j=1}|v_j|$ the $L_1$-norm. All vectors are column,  $\textbf{v}^t$ denotes the transposed of $\textbf{v}$ and $\frac{1}{\textbf{v}}=\pth{ \frac{1}{v_1}, \cdots, \frac{1}{v_p}}$.
For a strictly positive constant $c$, let also denote by $\textbf{v}^c$ the following vector $\pth{ \frac{1}{v_1^c}, \cdots, \frac{1}{v_p^c}}$. 
 Let $ \overset{\cal L} {\underset{n \rightarrow \infty}{\longrightarrow}}$, $ \overset{\eP} {\underset{n \rightarrow \infty}{\longrightarrow}}$, $ \overset{a.s.} {\underset{n \rightarrow \infty}{\longrightarrow}}$ represent convergence in distribution, in probability and almost sure, respectively, as $n \rightarrow \infty$. For a real $x$, $[x]$ is its integer part.   \\

The paper is organized as follows. In Section 2, we introduce and study especially the oracle properties of a general adaptive LASSO quantile estimator. In Section 3, the corresponding estimator in a change-point model is defined and its asymptotic behaviors (convergence rate and oracle properties) are studied. Section 4 proposes a consistent criterion to determine the breaking number. In Section 5, simulation results illustrate the performance of the proposed estimators. The adaptive LASSO quantile estimator is compared with other variable selection estimators.  All proofs are given in Section 6.

\section{Adaptive LASSO  quantile for regression model} 

In this section we propose and study a general  adaptive LASSO quantile estimator for the model (\ref{eq1}), estimator defined by
\begin{equation}
\label{eq4}
(\hat b^*_n, \hat \ef^*_n) \equiv \argmin_{(b,\ef)} \left(  \sum^n_{i=1} \rho_\tau (Y_i-b-\XX^t_i \ef) +\lambda_n \hat \eo^t_n | \ef| \right),
\end{equation}
where $\hat \eo_n=(\hat \omega_{n,1}, \cdots, \hat \omega_{n,p}) \equiv \frac{1}{| \hat \ef_n|^g}= \pth{ \frac{1}{| \hat \phi_{n,1}|^g}, \cdots,  \frac{1}{| \hat \phi_{n,p}|^g}}$, $ \hat \ef_n$ is the quantile estimator,  defined by (\ref{eq3}), and $g>0$ is a constant which will be  later specified. The components of $\hat \ef_n^*$ are $(\hat \phi^*_{n,1}, \cdots, \hat \phi^*_{n,p})$. The positive sequence $(\lambda_n)_n$, also called the tuning parameter, is a regularization parameter such that $\lambda_n \rightarrow \infty$ as $n \rightarrow \infty$. To the author's knowledge, for a quantile regression, three particular penalties of the  adaptive LASSO type have been previously proposed in other papers. In two papers, the case  $g=2$ is considered. First,    \cite{Zou:Yuan:08} consider that $p$, the number of regression parameters, is  fixed and after,  in  \cite{Belloni:Chernozhukov:11}, the case $p=p_n \rightarrow  \infty$ is studied. In the paper of \cite{Wu:Liu:09},   an adaptive LASSO penalty with the same form as in the relation (\ref{eq4}) is proposed, but under the assumption that the $\tau$th quantile $b^0_\tau$ is known and it is equal to 0. In the paper of \cite{Zheng:Gallagher:Kulasekera:13}, which is a particular case to that proposed here, always for  $b^0_\tau=0$,  the weight vector considered is $\min (\sqrt{n}, | \tilde \ef_n|^{-1})$, with $\tilde \ef_n$ a consistent estimator of $\efo$. However, their estimator has the advantage that it also applies when the number of regressors is very large, $p=O(\exp(n^c))$, with the constant $c \in (0,1)$.\\

Let us consider the deterministic design matrix $\mathbb{X}\equiv(X_{ij})_{\overset{1\leq i \leq n}{1 \leq j \leq p}}$,  with $\XX_i$ the $i$th line, corresponding to the observation $i$. \\
We give now two assumptions, denoted (A1), (A2), on the design and on the errors. \\ 
For the design  $\eeX$, we suppose that:\\
\textbf{(A1)} $n^{-1} \max_{1 \leq i \leq n} \XX_i^t \XX_i {\underset{n \rightarrow \infty}{\longrightarrow}}  0$ and 
$n^{-1} \sum^{n}_{i=1} \XX_i \XX_i^t {\underset{n \rightarrow \infty}{\longrightarrow}} \eU$, with $\eU=(v_{ij})_{1 \leq i,j \leq p}$ a positive definite matrix.\\
 For the errors $\varepsilon_i$ we suppose that they are  independent, identically distributed, with a continuous positive density $f$ in a neighborhood of $b^0_\tau$ and:\\
 \textbf{(A2)} For every $e \in int({\cal B})$, $u_0 \in \R$, $\uu \in \R^p$ we have 
 \[
 \lim_{n \rightarrow \infty} n^{-1} \sum^n_{i=1} \int^{u_0+\xx^t_i \uu}_{0} \sqrt{n}[F(e+v/\sqrt{n})-F(e)] dv = \frac{1}{2} f(e) (u_0,\uu^t) \left[ \begin{array}{cc}
 1 &  \textbf{0} \\
 0 & \eU
\end{array}  \right] (u_0,\uu^t)^t .
 \]
 Recall that  $F:{\cal B} \rightarrow [0,1]$ is the distribution function of $\varepsilon$. The $p$ regressors $X_{i1}, \cdots, X_{ip}$ are independent of the errors $\varepsilon_i$, for all $i=1, \cdots, n$. In fact, the design can be considered either deterministic or random, in which case we will consider the conditional expectation.  Obviously, the first  variable $X_1$ can be equal to 1, in which case the model contains an intercept. 
Remark that the assumption (A1) is standard when LASSO methods are used, see for example \cite{Zou:06}, \cite{Wu:Liu:09} and (A2) is classic for a quantile regression, see for example \cite{Koenker:05}. \\

Note that, the estimators defined by (\ref{eq3}) or (\ref{eq4}),  depend on the value of $\tau$ and on observation number $n$. For simplicity reasons, $\tau$ does not appear in the notation. 
The quantile estimators $(\hat b_n, \hat \ef_n)$ found by (\ref{eq3}), minimize the objective function $\sum^n_{i=1}[ \rho_\tau(\varepsilon_i-b-\XX_i^t(\ef-\efo))-\rho_\tau(\varepsilon_i-b^0_\tau)]$ and the adaptive LASSO quantile estimators $(\hat b^*_n, \hat \ef^*_n)$ given by (\ref{eq4}), minimize the objective function $\sum^n_{i=1}[ \rho_\tau(\varepsilon_i-b-\XX_i^t(\ef-\efo))-\rho_\tau(\varepsilon_i-b^0_\tau)]+\lambda_n \hat \eo^t_n [|\ef|-|\efo|]$. \\

To simplify the notations,  in what follows we denote $b^0_\tau$ by  $b^0$.\\

Then, let us first define two stochastic processes, necessary to study the objective functions and to analyze the behavior of the corresponding estimators obtained by quantile and by adaptive LASSO quantile methods, respectively.  
For $b \in {\cal B}$, $\ef \in \Gamma$ we   define, for each observation $i \in \{1, \cdots , n \}$, the following random process:
\begin{equation}
\label{Ri}
R^{(\tau)}_i (b,\ef;b^0,\efo) \equiv \rho_\tau(\varepsilon_i-b-\XX^t_i (\ef -\efo))-\rho_\tau(\varepsilon_i-b^0)
\end{equation}
and for two  observations $l$ and $k$ between 0 and $n$, with $l < k$, we define also the corresponding process, taking an adaptive LASSO penalty:
\[
R^{(\tau,\lambda)}_{i;(l,k)}(b,\ef;b^0,\efo) \equiv R^{(\tau)}_i (b,\ef;b^0,\efo) + \frac{\lambda_{(l,k)}}{k-l} \hat \eo^t_{(l,k)} (|\ef |-|\efo |), \qquad i=l+1, \cdots, k, 
\]
with  $\lambda_{(l;k)}$ the tuning parameter, dependent of position of  $l$ and $k$. The weight vector $\hat \eo_{(l;k)}=|\hat \ef_{(l;k)} |^{-g}$ is obtained considering  $\hat \ef_{(l;k)}$ the quantile estimator of $\ef$ calculated by (\ref{eq3}) on the samples $l+1, l+2, \cdots , k$: 
\[(\hat b_{(l;k)}, \hat \ef_{(l;k)})=\argmin_{(b,\ef)} \sum^k_{i=l+1} \rho_\tau(Y_i-b-\XX_i^t \ef).\]
 For $l=0$ and $k=n$, the tuning parameter $\lambda_{(0,n)}$ becomes $\lambda_n$  and $\hat \eo_{(0,n)}$, $(\hat b_{(0,n)}, \hat \ef_{(0,n)})$, $(\hat b^*_{(0,n)}, \hat \ef^*_{(0,n)})$ become $\hat \eo_n$, $(\hat b_n, \hat \ef_n)$, $(\hat b^*_n, \hat \ef^*_n)$ respectively.\\
\hh We note that the estimators obtained by  (\ref{eq3}) are in fact the  parameters that minimize in $(b,\ef)$ the random process $\sum^n_{i=1}R^{(\tau)}_i (b,\ef;b^0,\efo)  $  and the penalized estimators obtained  by relation (\ref{eq4}) are  the ones that minimize $\sum^n_{i=1} R^{(\tau,\lambda)}_{i;(0,n)}(b,\ef;b^0,\efo) $.\\

In the following Lemma we prove that, if in a model,  either the $\tau$th quantile or the  regression parameters are different to the true values, then the observation number for which $|b-b^0|+| \XX^t_i(\ef-\efo)|$ exceeds a threshold, is large. 
\begin{lemma} 
\label{Lemma 8Bai}
Under assumption (A1), if $b \neq b^0$ or $\ef \neq \efo$, then there exists $ \delta >0$ such that for the following set  
$ N_n \equiv Card \{ i \in \{1, \cdots, n\}; |b-b^0|+| \XX^t_i(\ef-\efo)| > \delta \}$ we have that there exists an  $\epsilon_0 >0$ such that   $N_n >n \epsilon_0$, for $n$ large enough.
\end{lemma}

In the following subsection we prove that, the adaptive LASSO estimator for the regression parameters satisfies the oracle properties: nonzero estimators are asymptotically normal and zero parameters are shrunk directly to 0 with a probability converging to 1 as $n$ tends to infinity. 

\subsection{Oracle properties}
First of all, let us formulate the  Karush-Kuhn-Tucker (KKT) optimality conditions, needed to prove the oracle properties. Note that for a real $x$, we use the notation $sgn(x)$ for the sign function $sgn(x) \equiv x/|x|$ when $x \neq 0$ and $sgn(0) =0$. Let  the set index 
\[\hat {\cal A}^*_n \equiv \{j \in \{1, \cdots , p \}; \hat \phi^*_{n,j} \neq 0   \}\]
 of nonzero components of the adaptive (LASSO)  quantile estimator of the regression parameter. 
\begin{proposition}
\label{KKT}
Under assumptions (A1), (A2), for $g >0$, if $(\lambda_n)$ is a sequence such that $\lambda_n \rightarrow \infty$, $n^{-1/2} \lambda_n \rightarrow 0$ and $n^{(g-1)/2} \lambda_n \rightarrow \infty$, as $n \rightarrow \infty$, then:\\
(i) for all $j \in \hat {\cal A}^*_n$ we have:
$
\tau \sum^n_{i=1} X_{ij} - \sum^n_{i=1} X_{ij} \e1_{Y_i < \XX^t_i \hat \ef^*_n}=n \hat \omega_{n,j}sgn(\hat \phi^*_{,j})$.\\
(ii) for all $j \not \in \hat {\cal A}^*_n$ we have $\left| \tau \sum^n_{i=1} X_{ij} - \sum^n_{i=1} X_{ij} \e1_{Y_i < \XX^t_i \hat \ef^*_n} \right| \leq n \hat \omega_{n,j}$.
\end{proposition}

In addition to the set $\hat {\cal A}^*_n$,  let us consider the set  
\[{\cal A}^0 \equiv \{ j \in \{1, \cdots, p \}; \phi^0_{,j} \neq 0\}  \]
 with the index of nonzero components of the true regression parameters. Throughout the paper, we denote by $\ef_{{\cal A}^0}$ the sub-vector of $\ef$ containing the corresponding components of ${\cal A}^0$. In  the same way, we consider   $\hat \ef^*_{\hat {\cal A}^*_n}$, $\hat \ef^*_{{\cal A}^0}$ subvectors of $\hat \ef^*_n$ with the index in the sets ${\hat {\cal A}^*_n}$, ${\cal A}^0$, respectively. \\
With these notations and with Proposition \ref{KKT}, we can now state the main result of this section on the oracle properties of the adaptive LASSO quantile estimator. 
We comment that even if the $\tau$th quantile is unknown, the assumptions are the same as in the known case of  \cite{Wu:Liu:09} where  $b^0$ was considered zero. 
In the paper of \cite{Zheng:Gallagher:Kulasekera:13} the particular case $\hat \omega_{n,j} = \min(\sqrt{n}, |\tilde \phi_{n,j}|^{-1}) $, with $\tilde \phi_{n,j}$ a consistent estimator of $\phi^0_{,j}$, is considered. In the present paper we consider $\hat \omega_{n,j} =|\hat \phi_{n,j}|^{-g} $, with $\hat \phi_{n,j}$ the quantile estimator.

\begin{theorem}
\label{Theorem 4.1 ZouYuan}
Under assumptions (A1), (A2), for a power $g >0$, if $(\lambda_n)$ is a sequence as in the Proposition \ref{KKT}, then the estimator $\hat \ef^*_n$ satisfies the oracle properties:\\
(i) $\sqrt{n} (\hat \ef^*_{{\cal A}^0}- \efo_{{\cal A}^0}) \overset{\cal L} {\underset{n \rightarrow \infty}{\longrightarrow}} {\cal N} \pth{ \textbf{0},   \frac{ \tau(1-\tau)}{f^2(b^0)}\eU^{-1}_{{\cal A}^0}} $, 
where   $\eU_{{\cal A}^0}$ contains the elements of the matrix $\eU$, defined in assumption  (A1), with the index in the set  ${\cal A}^0$: $\eU_{{\cal A}^0} \equiv (v_{ij})_{i,j \in {{\cal A}^0}}$.
\\
(ii) If moreover $n^{g/2-1} \lambda_n \rightarrow \infty$,  then $\eP[\hat {\cal A}^*_n={\cal A}^0] \rightarrow 1$,  as $n \rightarrow \infty$.
\end{theorem}

As a consequence of this Theorem, emphasize the fact that the convergence rate of the adaptive LASSO quantile estimators $\hat b^*_n$ and $\hat \ef^*_n$ has the order $n^{-1/2}$.\\
\hh The condition $n^{g/2-1} \lambda_n \rightarrow \infty$ on Theorem \ref{Theorem 4.1 ZouYuan}(ii) is stronger than that of Proposition \ref{KKT}, where $n^{(g-1)/2} \lambda_n \rightarrow \infty$, as $n \rightarrow \infty$, is supposed. The importance of this condition for  sparsity property will be studied in Section 5, by Monte Carlo simulations. Note that, the condition $n^{g/2-1} \lambda_n \rightarrow \infty$ on Theorem \ref{Theorem 4.1 ZouYuan}(ii) implies $g >1$.

\begin{remark}
In the particular case $b^0_\tau=0$, the Normal asymptotic law of  Theorem \ref{Theorem 4.1 ZouYuan}(i)   is the same that obtained by \cite{Wu:Liu:09} by the SCAD estimation method. If, moreover, $\tau=1/2$, we obtain the same law that \cite{Xu:Ying:10} using a LASSO-type penalty. 
\end{remark}
\begin{remark}
Compared with the adaptive LASSO estimator for a LS model (see \cite{Zou:06}), the fact that the conditions on $\varepsilon$ are weakened (no assumption on existence of the first two order moments and nor $\eE[\varepsilon]=0$), implies that in our case, in order to have the sparsity, the sequence  $(\lambda_n)$ and the  constant $g$ must verify the additional condition $n^{g/2-1} \lambda_n \rightarrow \infty$ beside  to claim (i).   
\end{remark}

\subsection{Behavior study for $R_i^{(\tau)}$ and $R_i^{(\tau,\lambda)}$}

The obtained  results  in this subsection on the random processes $R_i^{(\tau)}$ and $R_i^{(\tau,\lambda)}$ which  intervene in the corresponding objective functions, will be necessary to study  the behavior of the multiphase model.  
By the following  Proposition, we show that the process $R^{(\tau)}_i (b,\ef;b^0,\efo)$ has a positive expected value for any parameter $(b, \ef) \in {\cal B} \times \Gamma$.
\begin{proposition}
\label{proposition 2.1.SCAD}
For all parameters  $b \in {\cal B}$, $\ef \in \Gamma$, we have that $\eE[R^{(\tau)}_i(b,\ef;b^0,\efo)]\geq 0$.
\end{proposition}

 The penalty order, which is involved in the expression of $R_{i;(l;k)}^{(\tau,\lambda)}$, is obtained by the following result. It depends of the power $g$ considered for the weights $\hat \eo_{(l;k)}$.

 \begin{lemma}
 \label{lo}
 Under assumptions (A1), (A2), for all $0 \leq l <k \leq n$, if $\lambda_{(l;k)}=o(n^{1/2})$, then, $\lambda_{(l;k)} \| \hat \eo_{(l;k)}\|_2=o_{\eP}(n^{(1+g)/2})$.
 \end{lemma}

In order to study the multiphase model, the following result studies the objective function behavior between two observations sufficiently far apart. 

\begin{lemma}
\label{lemma 6aLASSO}
For $1 \leq l <k \leq n$ such that $k-l \rightarrow \infty $, as $n \rightarrow \infty$, under assumptions (A1) and (A2), if $\lambda_{(l;k)}=o(n^{1/2})$, for all $\alpha >1/2$, we have
$
\sup_{0 \leq l < k \leq n} \left| \inf_{b, \ef} \sum^{k}_{i=l} R_{i;(l;k)}^{(\tau,\lambda)}(b,\ef;b^0,\efo) \right|=\min(o_{\eP}(n^{(1+g)/2}),O_{P}(n^\alpha))$.
\end{lemma}

For notational simplicity, we will denote the sum of the processes $R_i^{(\tau)}$ on the observations $1, \cdots, n$ by: 
\[
{\cal R}_n^{(\tau)}(b,\ef;b^0,\efo) \equiv \sum^n_{i=1} R_i^{(\tau)}(b,\ef;b^0,\efo).
\] 
\hh Let us consider a positive deterministic sequence $(c_n)$, such that either it converges to 0, such that $ n c_n^2/ \log n \rightarrow \infty$ as $n \rightarrow \infty$, or it is constant $c_n=c$.  For a such sequence $(c_n)$, consider following sub-region of the set ${\cal B}\times \Gamma$: 
\[\Omega_n \equiv \{ (b, \ef ) \in {\cal B}\times \Gamma;  |b-b^0| \leq c_n, \| \ef-\efo\|_2 \leq c_n \}\] and $\Omega_n^c= \{(b, \ef ) \in {\cal B}\times \Gamma; \min(|b-b^0|, \| \ef -\ef^0\|_2 ) > c_n\} $ its complementary set. \\

The following Lemma shows that for all $n$ large enough, the supremum of the difference between  ${\cal R}_n^{(\tau)}$ and its expectation, when the parameters are in a neighborhood of the true values, converges to 0  in probability, with the rate $(n c^2_n)^{-1}$.
 \begin{lemma}
 \label{Lemma 4Bai}
Under assumptions (A1), (A2), then there exists a strictly positive constant $C_1>0$ such that for all $\epsilon >0$, there exists a $n_\epsilon \in \N$ such that, for $n \geq n_\epsilon$, following inequality holds
 \[
 \eP \cro{\sup_{(b,\ef) \in \Omega_n} \left|\frac{1}{n c^2_n} \cro{ {\cal R}_n^{(\tau)}(b,\ef;b^0,\efo)- \eE[{\cal R}_n^{(\tau)}(b,\ef;b^0,\efo)]}  \right| >\epsilon } \leq \exp (- \epsilon^2 n c_n^2 C_1).
 \]
 \end{lemma}

The following Lemma gives the behavior of the penalized objective function  $ \sum^n_{i=1}R^{(\tau,\lambda)}_{i;(0,n)} (b,\ef;b^0,\efo)$ on the complementary of the set $\Omega_n$. We obtain that the infimum of this process is strictly positive.  
 
 \begin{lemma}
 \label{lemma 7aLASSO}
Let us consider a positive deterministic sequence $(c_n)$ such that either it converges to 0, with $ n c_n^2/ \log n \rightarrow \infty$ as $n \rightarrow \infty$, or  $c_n=c$. If the tuning parameter sequence verifies $\lambda_n n^{-1} c^{-2}_n \rightarrow 0$, under assumptions (A1), (A2),  then we have that there exists $\epsilon_1 >0$ such that, with the probability 1:
 \[
 \liminf_{n \rightarrow \infty} \pth{ \inf_{(b,\ef) \in \Omega^c_n} \frac{1}{n c^2_n} \sum^n_{i=1} R^{(\tau,\lambda)}_{i;(0,n)} (b,\ef;b^0,\efo) } > \epsilon_1.
 \]
 \end{lemma}
 
 An example of a sequence $c_n$ and of tuning parameter $\lambda_n$ is $\lambda_n=n^{2/5}$ and $c^2_n=n^{-3/5} (\log n)^2$.  
  
\begin{remark}
\label{Remark 1}
Using Lemma \ref{lemma 7aLASSO}, by similar technique to one used in the paper of \cite{Ciuperca:13a}, for Lemma 3 and 4, we obtain  their equivalent. That is, if the data come from two different models, the adaptive LASSO quantile estimator is close to the parameter of the model from where most of the data came. 
\end{remark}

With these results, we can now consider a model with several phases. First, the number of phases is considered known, and afterward the number of changes will be assumed unknown. We will prove that the regression parameters of each phase and the change-points $l^0_{r-1}$, $l^0_r$ are estimated by consistent adaptive LASSO quantile estimators, with rate of convergence,  $(l_r^0-l^0_{r-1})^{-1/2}$,  $(l_r^0-l^0_{r-1})^{-1}$, respectively. A very interesting result is that oracle properties of the adaptive LASSO quantile estimators for the regression parameters  are not affected by the change-point estimation. 
  
\section{Adaptive LASSO  quantile for  multiphase model} 
Let us now consider a model with  $K+1$ phases, i.e. the model changes to the observations $l_1, \cdots , l_K$ with $1 < l_1 < l_2 < \cdots < l_K < n$. Initially, we suppose that the change number $K$ is known. If $K$ is unknown, which is the most frequent case in practice, we shall give in the following section a criterion to estimate the number $K$. \\
\hh The model with $K+1$ phases has the form:
\begin{equation}
\label{eq11}
Y_i = \XX^t_i \ef_1 \e1_{1 \leq i < l_1} + \XX^t_i \ef_2 \e1_{l_1 \leq i < l_2}+ \cdots +\XX^t_i \ef_{K+1} \e1_{l_K \leq i \leq n} +\varepsilon_i, \quad i=1, \cdots, n.
\end{equation}
\hh The parameters of model (\ref{eq11}) are $(b_1, \cdots, b_{K+1} )$  the  $\tau$th quantile of $\varepsilon$ on  each phase, $(\ef_1, \cdots, \ef_{K+1})$ the corresponding regression parameters and $(l_1, \cdots, l_K)$ the change-points (breaks). The true values of these parameters are  $(b^0_1, \cdots, b^0_{K+1} )$, $(\ef^0_1, \cdots, \ef^0_{K+1})$ and $(l^0_1, \cdots, l^0_K)$ respectively.\\
\hh Let us notice that with regard to the paper of  \cite{Ciuperca:13b} where  0 was always the $\tau$th quantile, here quantiles can change from one phase to another. In fact, in the present paper,  from one phase to the other, either the regression parameters  change, or the $\tau$th quantiles change, or both types of parameters changes simultaneously. To the author knowledge, this double possibility of change has not been addressed anywhere in the literature.  Moreover, in the paper of \cite{Ciuperca:13b}, the proposed penalty for a quantile multiphase model, with  the $\tau$th quantile  known, is SCAD (Smoothly Clipped Absolute Deviation), which produces difficulties of point of view numerical programming. In the particular case of a median multiphase model, \cite{Ciuperca:13b} considers a LASSO-type estimator in order to avoid the numerical  disadvantage generated by the SCAD method.    On the other hand, in view of simulations presented in Section 5, for median multiphase model, the results obtained by SCAD method are poorer than by a LASSO-type method. This justify the interest to consider an adaptive LASSO method for a quantile multiphase model. Moreover, the simulations presented in Section 5 will show the performance of the proposed method, especially in the case of not homoscedastic error or when the change occurs in the $\tau$th  quantile and not in the regression parameters.  \\

Concerning the  change-points $(l_1, \cdots , l_K)$, we suppose that a phase is long enough:\\
\textbf{(A3)} $l_{r+1} - l_r \geq n^a$, $a > 1/2$, for all $r=0, \cdots, K$, with $l_0=1$ and $l_{K+1}=n$.\\

For fixed change-points $(l_1, \cdots , l_K)$, the objective function is minimized with respect to the $\tau$th quantiles and regression parameters of the $K+1$ phases. Let us denote the objective function value by 
\begin{equation}
\label{def_S*}
S^*(l_1, \cdots, l_K)\equiv \inf_{
\begin{array}{c}
(\ef_1, \cdots , \ef_{K+1})\\ 
(b_1, \cdots, b_{K+1})
\end{array}} \sum^{K+1}_{r=1} [ \sum^{l_r}_{i=l_{r-1}+1} \rho_\tau(Y_i-b_r-\XX^t_i \ef_r) + \lambda_{(l_{r-1};l_r)} \hat \eo^t_{(l_{r-1};l_r)} | \ef_r|].
\end{equation}
The tuning parameters $\lambda_{(l_{r-1};l_r)}$ vary from a phase to the  other one with the interval length $l_r -l_{r-1}$. The weight vector is $\hat \eo^t_{(l_{r-1};l_r)}= |\hat \ef_{(l_{r-1}; l_r)} |^{-g}$, with $\hat \ef_{(l_{r-1}; l_r)}$ the quantile estimator of $\ef_r$ calculated by quantile method, on the observations $l_{r-1}+1, \cdots, l_r$. \\
\hh Then the adaptive LASSO quantile estimators for the change-points are the minimizers of the function $S^*$:
\[(\hat l^*_1, \cdots,\hat  l^*_K)\equiv \argmin_{(l_1, \cdots, l_K)} S^*(l_1, \cdots, l_K). 
\]
The adaptive LASSO quantile estimator for  the regression parameters is $\hat \ef^*_{(\hat l^*_{r-1};\hat l^*_r)}$ and for the $\tau$th quantile is  $\hat b^*_{(\hat l^*_{r-1}; \hat l^*_r)}$, for each $r=1, \cdots, K+1$:\\
\[
\begin{array}{l}
((\hat b^*_{(\hat l^*_{0}; \hat l^*_1)}, \hat \ef^*_{(\hat l^*_{0};\hat l^*_1)}), \cdots , (\hat b^*_{(\hat l^*_{K}; \hat l^*_{K+1})}, \hat \ef^*_{(\hat l^*_{K};\hat l^*_{K+1})}) )=\\
 \argmin_{
\begin{array}{c}
(\ef_1, \cdots , \ef_{K+1})\\ 
(b_1, \cdots, b_{K+1})
\end{array}} \sum^{K+1}_{r=1} [ \sum^{\hat l^*_r}_{i=\hat l^*_{r-1}+1} \rho_\tau(Y_i-b_r-\XX^t_i \ef_r) + \lambda_{(\hat l^*_{r-1};\hat l^*_r)} \hat \eo^t_{(\hat l^*_{r-1};\hat l^*_r)} | \ef_r|].
\end{array}
\]

In the next result we show that if in a phase we take in the place of the true parameters   those of the nearby phase, then,  we obtain, with a  probability close to  1, different values  for the objective function (without penalty).  
  
\begin{lemma} 
\label{Lemma 9Bai}
Under assumptions (A1), (A2), for $r=1, \cdots, K$, if $l_r$ is such that $l_r < l^0_r$, $l_r-l^0_r=O(1)$, we have that there exists two strictly  positive constants 
$\eta $, $ C_1 $, such that
\[
\eP \cro{\sum^{l^0_r}_{i=l_r+1}[\rho_\tau(\varepsilon_i-b^0_{r+1}-\XX^t_i(\ef^0_{r+1}-\ef^0_r)) - \rho_\tau(\varepsilon_i-b^0_r) ] \geq \eta(l^0_r - l_r) } \geq 1-\exp(-C_1(l^0_r-l_r)).
\]
\end{lemma}

  %%%%%%%%%%%%%%%%%%%%%%%%%%%%%%%%%%%%%%%ù
%\begin{figure}[h!]
%\begin{minipage}[b]{0.50\linewidth}
%\includegraphics[width=10cm,height=1.5cm,angle=0]{fig0.pdf} 
%\includegraphics[width=10cm,height=15.5cm,angle=0]{nlin2_laplace.eps}
%\caption{\it Positionement change-points.}
%\label{Figure 1}
%\end{minipage}
%\end{figure}

%%%%%%%%%%%%%%%%%%%%%%%%%%%%%%%%%%%%
Let us now formulate the equivalent to assumption (A1) in the case of a model with $(K+1)$ phases, when the change-points are sufficiently far apart, i.e. $l_r-l_{r-1}$ converges to infinity as $n \rightarrow \infty$:\\

\textbf{(A1bis)} $(l_r -l_{r-1})^{-1} \max_{l_{r-1} < i \leq l_r} \XX^t_i \XX_i \rightarrow 0$ for $n\rightarrow \infty $. We suppose that for each phase we have that the matrix $(l_r-l_{r-1})^{-1} \sum^{l_r}_{i=l_{r-1}+1} \XX_i \XX^t_i$ converges to  $ \eU_r$, as $n \rightarrow \infty$, with $\eU_r$ a positive definite matrix. Let us denote by $\eU^0_r$ the limiting matrix for the true change-points $l^0_r$, $r=1, \cdots, K$. We also denote by $v^0_{r,kj}$ the $(k,j)$th component of matrix $\eU^0_r$.\\ 

Following result gives the convergence rate of the adaptive LASSO quantile estimator of the change-points.  The proof is quite technical, that is why we will divide it into three parts. You can not directly show that the distance between the estimator and the true value of the corresponding change-point is finished. This is why we first show that this distance is less than $[n^{1/2}]$, afterward  less than  $[n^{1/4}]$ and at the end that it is bounded.
 \begin{theorem}
 \label{Theorem 5aLASSO}
 Under  assumptions (A1bis), (A2), (A3),  if the tuning parameter $(\lambda_{(l_{r-1},l_r)})_{1 \leq r \leq K+1}$   is a sequence, depending on $n$, converging to zero, such that  $(l_r-l_{r-1})^{1/2}\lambda_{(l_{r-1},l_r)} \rightarrow \infty$ and if $(c_n)$ is another  deterministic sequence $(c_n)$, such that   $c_n \rightarrow 0$, $n c^2_n/ \log n\rightarrow \infty$ and $\lambda_n c^{-2}_n  \rightarrow 0$, as $n \rightarrow \infty$, then  we have 
 $\hat l^*_r- l^0_r=O_{\eP}(1)$ for each $r=1, \cdots, K$.
 \end{theorem}
 
 If instead of the parameters $l_1, \cdots, l_K$ we consider the reparametrization $\theta_r=l_r/n \in (0,1)$, with $0 < \theta_1 < \theta_2 < \cdots \theta_K  < 1$, the adaptive LASSO  quantile estimators are   $\hat \theta^*_r =\hat l^*_r/n$, which have a convergence rate of order $n^{-1}$. This is the classical rate convergence for the change-point estimators, when the number of regressors is not dependent on $n$. See for example the paper of \cite{Bai:98} for the  LAD estimators, or the paper of  \cite{Ciuperca:11} for penalized  LAD estimators, but not with adaptive type or LASSO penalty. For the  LS method, without penalty, the reader can see \cite{Bai:Perron:98}. \\
 \hh Concerning  the estimators  $\hat \ef^*_{(\hat l^*_{r-1};\hat l^*_r)}$, in view of  the results obtained in Section 2 and of Theorem \ref{Theorem 5aLASSO}, we have that its convergence rate is of order  $(l^0_r-l^0_{r-1})^{1/2}$.\\
 
 For $r=1, \cdots , K+1$ and $j=1, \cdots , p $, denote by:\\
 $ \phi_{r,j}$  the  $j$th component of the true value  $\ef^0_r$, \\
 $ \hat \phi^*_{(l^0_{r-1};l^0_r),j}$  the  $j$th component of the adaptive LASSO estimator $\hat \ef^*_{(l^0_{r-1};l^0_r)}$,\\ 
 $ \hat \phi^*_{(l^*_{r-1};l^*_r),j} $ the $j$th component of the adaptive LASSO estimator  $\hat \ef^*_{(l^*_{r-1};l^*_r)}$.\\ 
  
The following Theorem shows that the oracle properties are preserved in a multiphase model. For each two consecutive true change-points $l^0_{r-1}$, $l^0_r$ consider the set with the index of nonzero components of the true regression parameters 
\[{\cal A}^0_r \equiv \{ j \in \{1, \cdots, p\};  \phi^0_{r,j}\neq 0 \}\]
 and following set when only the regression parameters are estimated  
 \[\hat {\cal A}^0_{n,r} \equiv \{ j \in \{1, \cdots, p\}; \hat \phi^*_{(l^0_{r-1},l^0_r),j}\neq 0 \}.\]
  Consider also the similar index set corresponding to the adaptive LASSO quantile  estimators $\hat l^*_{r-1}$, $\hat l^*_r$ of the change-points: 
  \[\hat {\cal A}^*_{n,r} \equiv \{ j \in \{1, \cdots, p\}; \hat \phi^*_{(\hat l^*_{r-1},\hat l^*_r),j}\neq 0 \}.\]
In the following theorem, the matrix $\eU^0_{ {\cal A}^0_r}$ contains the  elements of matrix $\eU_r$, the vectors $\pth{\hat \ef^*_{(\hat l^*_{r-1};\hat l^*_r)}-\ef^0_r}_{ {\cal A}^0_r}$,  $\pth{\hat \ef^*_{(\hat l^*_{r-1};\hat l^*_r)}-\ef^0_r}_{ {\cal A}^0_r}$ are subvectors of $\pth{\hat \ef^*_{(\hat l^*_{r-1};\hat l^*_r)}-\ef^0_r}$, $\pth{\hat \ef^*_{(\hat l^*_{r-1};\hat l^*_r)}-\ef^0_r}$, respectively, all    with the index in the set $ {\cal A}^0_r$.
\begin{theorem}
\label{Theorem oracle}
 Under  assumptions (A1bis), (A2), (A3), for a power  $g >0$  the tuning parameter sequence $(\lambda_{(l_{r-1},l_r)})_{1 \leq r \leq K+1}$ on each interval $(l_{r-1},l_r)$ as in Theorem \ref{Theorem 5aLASSO} and also $(l_r-l_{r-1})^{(g-1)/2} \lambda_{(l_{r-1}, l_r)} \rightarrow \infty$, as $n \rightarrow \infty$, we have \\
(i) $(\hat l^*_r-\hat l^*_{r-1})^{1/2} \pth{\hat \ef^*_{(\hat l^*_{r-1};\hat l^*_r)}-\ef^0_r}_{ {\cal A}^0_r}=(l^0_r-l^0_{r-1})^{1/2}\pth{\hat \ef^*_{(\hat l^*_{r-1};\hat l^*_r)}-\ef^0_r}_{ {\cal A}^0_r}(1+o_{\eP}(1)) \overset{{\cal L}} {\underset{n \rightarrow \infty}{\longrightarrow}}$ $ {\cal N} \pth{\oo, (\tau(1-\tau))f^{-2}(b^0_r)(\eU^0_{{\cal A}^0_r})^{-1} }$.  \\
(ii) If the tuning parameter $\lambda_{(l_{r-1}, l_r)}$ satisfies more $(l_r-l_{r-1})^{g/2-1} \lambda_{(l_{r-1}, l_r)} \rightarrow \infty $, as $n \rightarrow \infty$, then, for every $r=1, \cdots, K$, we have $\lim_{n \rightarrow \infty } \eP \cro{\hat {\cal A}^0_{n,r} = \hat {\cal A}^*_{n,r} = {\cal A}^0_r}=1$.
\end{theorem} 

  %%%%%%%%%%%%%%%%%%%%%%%%%%%%%%%%
 % \begin{figure}[h!]
%\begin{minipage}[b]{0.48\linewidth}
%\includegraphics[scale=0.45]{fig1.pdf} 
%\includegraphics[width=10cm,height=15.5cm,angle=0]{nlin2_laplace.eps}
%\caption{\it Figure 1.}
%\label{Figure 1}
%\end{minipage}
%\begin{minipage}[b]{0.53\linewidth}
%\includegraphics[scale=0.45]{fig2.pdf}
%\includegraphics[width=10cm,height=15.5cm,angle=0]{nlin3_laplace.eps}
%\caption{\it Figure 2.}
%\label{Figure 2}
%\end{minipage}
%\end{figure}
%%%%%%%%%%%%%%%%%%%%%%%%%%%%%%%%%%%%%%%%

   \section{Selection criterion of the number of phases}
   
 Let us now give a criterion to estimate the true change-points number, noted  $K^0$,  of the model (\ref{eq11}).\\
 First give some notations and an additional assumption.  \\
 
  Concerning the distribution of error  $(\varepsilon_i)$, we suppose that\\
   \textbf{(A4)} $0 < \eE[\rho_\tau(\varepsilon-b_r^0)]$ and $ \eE[\rho^2_\tau(\varepsilon-b_r^0)]  < \infty$ for all $r= 1, \cdots, K^0$.\\
   
   If  $b^0_r=0$, then the first part of the assumption (A4) amounts to imposing that the mean of  $|\varepsilon|$ is bounded and  strictly positive. In fact, since $|\rho_\tau(\varepsilon-b^0_r)-\rho_\tau(\varepsilon)| \leq  |b^0_r |$ with probability 1, and since  $b^0_r$ belongs to a compact set, the assumption  (A4) implies that $\eE[\rho_\tau(\varepsilon)] < \infty$.  \\
The condition $ \eE[\rho^2_\tau(\varepsilon-b_r^0)]  < \infty$ is necessary to define a consistent criterion  in case when the error distribution changes at the  observation $l^0_r$, for $r=1, \cdots, K^0$.\\

For any number $K$ of changes, we calculate the sum  $S^*(l_1, \cdots, l_K)$, defined by the relation (\ref{def_S*}),  and the corresponding change-point estimators thereby: 
\[(\hat l^*_{1,K}, \cdots , \hat l^*_{K,K}) \equiv \argmin_{(l_1, \cdots, l_K)} S^*(l_1, \cdots, l_K).\]
Let us consider the objective function divided by the observation number
$\hat s^*_K \equiv n^{-1} S^*(\hat l^*_{1,K}, \cdots , \hat l^*_{K,K})$.
 In order to find the estimator of $K$, let us consider following criterion
\begin{equation}
\label{eq18}
B(K) \equiv n \log \hat s^*_K+G(K,p_K) B_n,
\end{equation}
where $(B_n)$ is a deterministic sequence converging  to infinity such that  $B_n n^{-a} \rightarrow 0$, $B_n n^{-1/2} \rightarrow \infty$ as $n \rightarrow \infty$. The constant  $a$ is that of the  supposition (A3). The penalty $G(K,p)$ is a function such that $G(K_1,p) \leq G(K_2,p)$, for all $K_1 \leq K_2$ and optionally depending on the number $p$ of  parameters to be estimated. \\
\hh By the proof of  Theorem \ref{Theorem 4aLASSO}, we have that $\hat s^*_{K^0}$ converges in probability to $\sum^{K^0+1}_{r=1} \eE[\rho_\tau(\varepsilon +b^0_r) ]$. Then, in order that the proposed criterion is well defined, it is necessary to impose the condition that each $\eE[\rho_\tau(\varepsilon +b^0_r)] $ is strictly positive and bounded.   \\
\hh We consider as an estimator for $K^0$, the number of change-points that minimizes the criterion $B(K)$, so
\begin{equation}
\label{eq19}
\hat K^*_n \equiv \argmin_K \pth{n \log \hat s^*_K +G(K,p) B_n }. 
\end{equation}
\hh This type of criterion for choosing the change-point number  was introduced, as a  Schwarz criterion, by \cite{Yao:88} for a constant model, scalar, in each phase, with $G(K,p)=K$. It was after considered, for an without penalty median model, by \cite{Bai:98}.  Other information criterion are used to detect the change number in the papers  \cite{Osorio:Galea:05}, \cite{Wu:08}, \cite{Nosek:10}.  In the paper of \cite{Liu:Zou:Zhang:08}, the empirical likelihood test was considered to detect a single change against no-change in a linear regression. \cite{Qu:Perron:07} consider likelihood ratio type statistic to test the null hypothesis of $K$ changes, against the alternative hypothesis of $(K+1)$ changes under the assumptions that the errors have mean 0 and bounded variance. But this approach type has the disadvantage that it must perform successive test to find the true number $K^0$ of the change-points. The only case when two changes are tested agains no changes is the particular case of the epidemic charge model. See for example the paper of \cite{Ning:Pailden:Gupta:12} for the epidemic change in a constant model. \\

We prove now  that $\hat K^*_n$ obtained by the relation  (\ref{eq19}) is a weakly convergent estimator for $K^0$. 

\begin{theorem}
\label{Theorem 4aLASSO}
Under assumptions (A1bis), (A2)-(A4), if the deterministic  sequence  $(B_n)$ converging to infinity  is such that   $n^{-a}B_n \rightarrow 0$, $n^{-1/2} B_n  \rightarrow \infty$, as $n \rightarrow \infty$, then we have that $\hat K^*_n \overset{\eP} {\underset{n \rightarrow \infty}{\longrightarrow}} K^0$.
\end{theorem}

\begin{remark}
In view of the proof of  Theorem \ref{Theorem 4aLASSO}, for a fixed index $\tau$,  when the error quantile does not change from one phase to the other,  then for assumption (A4) we require only $0 < \eE[\rho_\tau(\varepsilon-b_r^0)] < \infty$.
\end{remark}

\section{Simulations}
All simulations were performed using the R language. The program codes are available from the author. First we compare the results of the proposed method with other existing in the literature for a model with a single phase. After, the best three methods are studied for a multiphase model.
\subsection{Models with a single phase}
The samples were generated from the following model with a single phase: $Y_i=\XX_i^t \efo +\varepsilon_i$, $i=1, \cdots , n$, with $\efo=(1,0,4,0,-3,5,6,0,-1,0)$, $\XX=(X_1, \cdots, X_{10})$,  $X_3 \sim {\cal N}(2,1)$, $X_4 \sim {\cal N}(-1,1)$, $X_5 \sim {\cal N}(1,1)$ and $X_j \sim {\cal N}(0,1)$ for $j \in \{1,2,6,7,8,9,10 \}$. For the errors $\varepsilon$, three distributions were first considered: standard Normal ${\cal N}(0,1)$, exponential ${\cal E}xp(-4.5,1)$ with the density function $\exp (-(x+4.5)) \e1_{x > -4.5}$, and Cauchy ${\cal C}(0,1)$. The sample size is $n=200$. \\
\hh The percentage of zero coefficients correctly estimated to zero (true 0) and the percentage of nonzero coefficients estimated to zero (false 0) are computed by four methods: least squares model with adaptive LASSO penalty, median model with LASSO-type penalty, quantile model with adaptive LASSO and SCAD penalties. \\
 Recall (see the paper of \cite{Zou:06}) that the adaptive LASSO estimators for the regression parameters in a LS model, are the minimizers of the following objective function 
 \[\sum^n_{i=1} (Y_i-\XX_i \ef)^2 +\lambda_n \hat \varpi_n |\ef|,\]
  with the adaptive weight p-vector $\hat \varpi_n$  considered here that $|\hat \ef^{LS}_n|^{-\chi }$. Precise that  $\hat \ef^{LS}_n$ is the LS estimator of $\ef$. The adaptive LASSO estimator  for  LS model has the sparsity property if $\lambda_n \rightarrow \infty $, $n^{-1/2} \lambda_n \rightarrow 0$ and $n^{(\chi-1)/2} \lambda_n \rightarrow \infty$, as $n \rightarrow \infty$.  We consider then the tuning parameter  $\lambda_n=n^{2/5}$ and the power $\chi=9/40$. \\
 The SCAD estimator for a quantile model (see \cite{Wu:Liu:09}) is the minimizer of the following objective function
 \[\sum^n_{i=1} [ \rho_\tau(Y_i-\XX^t \ef)+\sum^p_{j=1} p_{\lambda_n} (| \phi_{,j}|)],\] 
 with penalty $p_{\lambda_n}(|\phi_{,j}|)$  defined by its first derivative
\[p'_{\lambda_n}(|\phi_{,j}|) = \lambda_n \{\e1_{|\phi_{,j}| \leq \lambda_n} +\frac{(a_1 \lambda_n -|\phi_{,j}|)_{+}}{(a-1) \lambda_n} \e1_{|\phi_{,j}| > \lambda_n} \},\]
  for all $j=1, \cdots, p$,
with $\lambda_n >0$, $a_1>2$ deterministic tuning parameters. We take here, for the SCAD method, the tuning parameters $a_1=5$ and $\lambda_n= 1/| \hat \ef^{QLASSO}_n|$. The estimator $\hat \ef^{QLASSO}_n$ of $\ef$ is the minimizer of the objective function $\sum^n_{i=1} [\rho_\tau(Y_i-\XX^t \ef)+ \el_n |\ef|$, with $\el_n=\log n \cdot \uuu_p$. \\
The LASSO-type estimator for a median model (see \cite{Xu:Ying:10}), is the minimizer of the following objective function 
\[\sum^n_{i=1} | Y_i- \XX^t_i \ef|+\el_n^t | \ef |,\]
 with the tuning parameter $\el_n$ a random $p$-vector. We take here $\el_n=n^{2/5} \frac{1}{|\hat \ef^{QLASSO}_n|}$. \\
\hh For adaptive LASSO quantile method, given by relation (\ref{eq4}),  we consider two values for the power $g$: a value greater than 1 and another smaller than 1. The tuning parameter is $\lambda_n=n^{2/5}$. \\
\hh The results obtained by these four methods are presented in Tables \ref{Tabl1}-\ref{Tabl3} for 1000 Monte Carlo replications. For some distributions or some quantile index, there is numerical problems (the function \textit{rq} of the package \textit{quantreg} of R language does not respond) for the SCAD method. Then, in Tables \ref{Tabl1}-\ref{Tabl3}, this is symbolised by "???". The numerical problems of the SCAD method have been also identified by \cite{Xu:Ying:10} who proposed the LASSO-type penalty for median regression.

\subsubsection{Sparsity property} 

For the adaptive LASSO quantile method proposed in this paper, considering  three error distributions, we deduce  from Tables \ref{Tabl1}-\ref{Tabl3} that the sparsity is not satisfied when the power $g$ of the adaptive weight is smaller than 1. This is in concordance with  the condition imposed in statement (ii) of the Theorem \ref{Theorem 4.1 ZouYuan}: $n^{g/2-1} \lambda_n \rightarrow \infty$. In all of the above simulations, we have the following conclusions: for $g>1$, the performance of  the adaptive LASSO quantile estimations are always better than the SCAD, median (LAD) model with LASSO-type penalty  and adaptive LASSO (for LS model) estimators,  for the heavy-tailed errors.

 \begin{table}
{\scriptsize
\caption{\footnotesize Model with a single phase.  Percentage  true 0 and of false 0 by LS+adaptiveLASSO, QUANTILE+adaptiveLASSO, QUANTILE+SCAD, LAD+LASSOtype methods for $n=200$,  $\varepsilon_i \sim {\cal N}(0,1)$.}
\begin{center}
\begin{tabular}{|cc|c|c|c|c|c|} \hline 
$\tau \downarrow$ & Method $\rightarrow$& LS+aLASSO & QUANT+aLASSO & QUANT+aLASSO  & QUANT+SCAD & LAD+LASSOtype \\ 
 & parameters $\rightarrow$ & $\chi=9/40$ & $g_1=12.25/10$ & $g_2=9/40$  &  &\\ \hline
0.15   &$\%$ of trues 0 &  1 & 1 & 0.77 & 0.47 & 0.99\\
   &$\%$ of   false 0 & 0.005& 0 & 0 & 0 & 0 \\ \hline
   0.50   &$\%$ of trues 0 &  1 & 1 & 0.64 & 0.55 & 0.99 \\
   &$\%$ of   false 0 & 0.01 & 0 & 0 & 0 & 0  \\ \hline
   0.95   &$\%$ of trues 0 &  1 & 1 & 0.93 & ??? & 0.99 \\
   &$\%$ of   false 0 & 0.01& 0 & 0 & ??? & 0 \\ \hline
\end{tabular} 
\end{center}
\label{Tabl1} 
}
\end{table}

\begin{figure}[h!]
   \begin{minipage}[b]{0.55\linewidth}
\centering \includegraphics[scale=0.40]{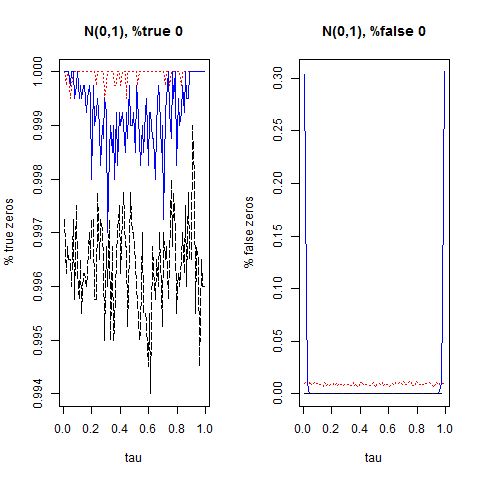}
 % laplace_lin.eps: 1179666x1179666 pixel, 300dpi, 9987.84x9987.84 cm, bb=   67   212   552   589
\caption{ \it  Percentage of  true 0 and of false 0 by LS+aLASSO, QUANTILE+aLASSO, LAD+LASSOtype methods, $\varepsilon_i \sim  {\cal N}(0,1)$.}
\label{Figure 4}  
   \end{minipage}\hfill
\begin{minipage}[b]{0.55\linewidth}
 \centering \includegraphics[scale=0.40]{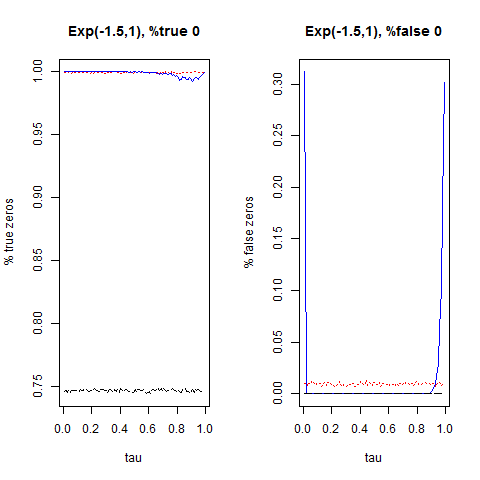}
 % normal_lin.eps: 6298172x0 pixel, 300dpi, 53324.52x0.00 cm, bb=   78   212   546   589

\caption{\it  Percentage of true 0 and of false 0 by LS+aLASSO, QUANTILE+aLASSO, LAD+LASSOtype methods, $\varepsilon_i \sim {\cal E}xp(-4.5,1)$.}
\label{Figure 5}
\end{minipage}
\end{figure}

\begin{figure}[h!]
   \begin{minipage}[b]{0.55\linewidth}
\centering \includegraphics[scale=0.40]{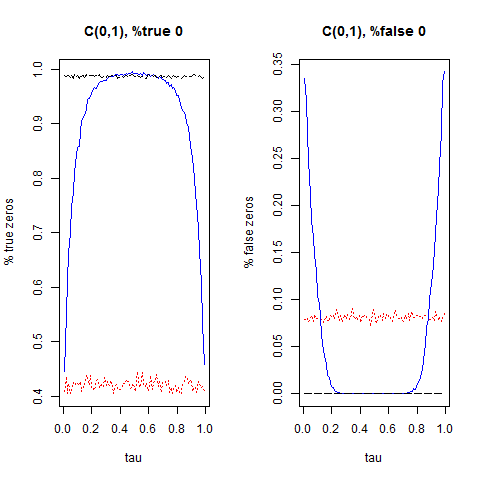}
 % laplace_lin.eps: 1179666x1179666 pixel, 300dpi, 9987.84x9987.84 cm, bb=   67   212   552   589
\caption{ \it  Percentage of  true 0 and of false 0 by LS+aLASSO, QUANTILE+aLASSO, LAD+LASSOtype methods, $\varepsilon_i \sim {\cal C}(0,1)$.}
\label{Figure 6}  
   \end{minipage}\hfill
%\begin{minipage}[b]{0.55\linewidth}
%\end{minipage}
\begin{minipage}[b]{0.55\linewidth}
 \centering \includegraphics[scale=0.38]{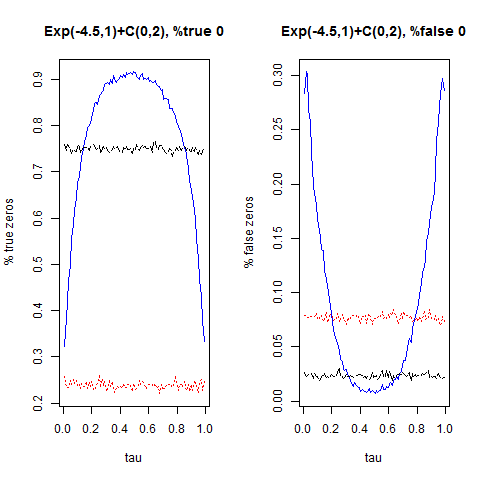}
 % normal_lin.eps: 6298172x0 pixel, 300dpi, 53324.52x0.00 cm, bb=   78   212   546   589

\caption{\it  Percentage of true 0 and of false 0 by LS+aLASSO, QUANTILE+aLASSO, LAD+LASSOtype methods, $\varepsilon_i \sim {\cal E}xp(-4.5,1)+C(0,2)$.}
\label{Figure 7}
\end{minipage}
\end{figure}

\begin{table}
{\scriptsize
\caption{\footnotesize Model with a single phase.  Percentage  true 0 and of false 0 by adaptive LASSO least squares, adaptive LASSO quantile, QUANTILE+SCAD, LAD+LASSOtype  methods for $n=200$,  $\varepsilon_i \sim {\cal E}xp(-4.5,1)$.}
\begin{center}
\begin{tabular}{|cc|c|c|c|c|c|} \hline 
$\tau \downarrow$ & Method $\rightarrow$& LS+aLASSO & QUANT+aLASSO & QUANT+aLASSO  & QUANT+SCAD & LAD+LASSOtype \\ 
 & parameters $\rightarrow$ & $\chi=9/40$ & $g_1=12.25/10$ & $g_2=9/40$  & & \\ \hline
0.15   &$\%$ of trues 0 &  0.99 & 1 & 0.88 & 0.28 & 0.66 \\
   &$\%$ of   false 0 & 0.01& 0 & 0 & 0 &  0 \\ \hline
   0.50   &$\%$ of trues 0 &  1 & 1 & 0.71 & 0.28 & 0.67 \\
   &$\%$ of   false 0 & 0.01& 0 & 0 & 0  & 0 \\ \hline
   0.95   &$\%$ of trues 0 &  1 & 0.99 & 0.90 & ??? & 0.67 \\
   &$\%$ of   false 0 & 0.01& 0.02 & 0.01 & ???  & 0 \\ \hline
\end{tabular} 
\end{center}
\label{Tabl2} 
}
\end{table}
\begin{table}
{\scriptsize
\caption{\footnotesize  Model with a single phase.  Percentage  true 0 and of false 0 by LS+adaptiveLASSO, QUANTILE+adaptiveLASSO, QUANTILE+SCAD,  LAD+LASSOtype methods for $n=200$, $\varepsilon_i \sim C(0,1)$.}
\begin{center}
\begin{tabular}{|cc|c|c|c|c|c|} \hline 
$\tau \downarrow$ & Method $\rightarrow$& LS+aLASSO & QUANT+aLASSO & QUANT+aLASSO  & QUANT+SCAD & LAD+LASSOtype \\ 
 & parameters $\rightarrow$ & $\chi=9/40$ & $g_1=12.25/10$ & $g_2=9/40$  & &  \\ \hline
0.15   &$\%$ of trues 0 &  0.43 & 0.93 & 0.69 & ??? & 0.99 \\
   &$\%$ of   false 0 & 0.08 & 0.03 & 0.025 & ???  & 0 \\ \hline
   0.50   &$\%$ of trues 0 &  0.44 & 0.99 & 0.64 & ??? & 0.98\\
   &$\%$ of   false 0 & 0.08& 0 & 0 & ???  & 0   \\ \hline
   0.95   &$\%$ of trues 0 &  0.43 & 0.66 & 0.77 & ???  & 0.98  \\
   &$\%$ of   false 0 & 0.08& 0.24 & 0.28 & ??? &  0 \\ \hline
\end{tabular} 
\end{center}
\label{Tabl3} 
}
\end{table}
 In view of the obtained results, presented in Tables \ref{Tabl1}-\ref{Tabl3}, the penalty SCAD and the parameter $g <1$ for the adaptive LASSO quantile methods are abandoned.\\
\hh Based on 1000 Monte-Carlo replications, in  Figures \ref{Figure 4}, \ref{Figure 5}, \ref{Figure 6}, \ref{Figure 7} we represent the graph of the percentages of true 0 and of false 0 obtained by three estimation methods:
 \begin{itemize}
 \item for LS model with adaptive LASSO penalty (dotted line),
 \item for quantile model with adaptive LASSO penalty (solid line), 
 \item for median model with LASSO-type penalty (long dash line),
\end{itemize}    
each for errors with four possible distributions: Normal, Exponential,  Cauchy and ${\cal E}xp(-4.5,1)+{\cal C}(0,2)$. \\
For Gaussian errors, these three methods give very satisfactory results (see Figure \ref{Figure 4}). For Exponential errors ${\cal E}xp(-4.5,1)$, the LASSO-type method does not well identify the true zeros (see Figure \ref{Figure 5}) since this method is build for median regression. For Cauchy errors ${\cal C}(0,1)$ (then moments of errors don't exist) the adaptive LASSO method for LS model provides estimates that don't have the sparsity property (see Figure \ref{Figure 6}). The superiority of the adaptive LASSO quantile method is considerably higher when errors are a sum of Exponential and Cauchy laws (see Figure \ref{Figure 7}).
 
We note that, in Figures \ref{Figure 4}-\ref{Figure 7}, for the LS model with adaptive LASSO penalty (dotted line), the quantile index $\tau$ is not useful, the represented values of the true or false zeros are in fact Monte Carlo replications.

\subsubsection{Conclusion}

In conclusion, for an uniphase model, adaptive LASSO quantile method provides very satisfactory results of sparsity for each distribution error, tacking as quantile index $\tau \in [0.4; 0.6]$.\\
\hh  For the sparsity property, the power $g$ of the relation (\ref{eq4}), should be greater than 1. Condition on $g$ in concordance with the condition imposed in statement (ii) of Theorem \ref{Theorem 4.1 ZouYuan}.\\
\hh  The quantile model with SCAD penalty gives bad sparsity results, and moreover,  it has numerical problems. \\
\hh Concerning the sparsity, the adaptive LASSO quantile model, the adaptive LASSO method for LS model and LASSO-type method for median model, give very satisfactory results for Gaussian errors. The adaptive LASSO quantile method stands out to be the best method, in terms of variable selection, when moments of errors don't exist or when median of errors is different to zero.
\subsection{Multiphase models} 
In view of  the results for uniphase model, we will consider in this subsection only three estimation methods: adaptive LASSO quantile, adaptive LASSO for LS model and LASSO-type for median model.
For each of the three methods, we have considered the same powers, tuning parameters as is the previous sub-section, with the only difference that instead of $n$ we take $l_r-l_{r-1}$, for $r=1, \cdots , K+1$.
 The quantile index $\tau$ for quantile model with  adaptive LASSO penalty is considered 0.55.\\
 
\subsubsection{Fixed phase number}
We consider that the change number  is known and it is equal to two (three phases). The true change-points are in $l^0_1=30$, $l^0_2=100$ for $n=200$ observations. In Tables \ref{Tabl5} and \ref{Tabl5bis}, we present simulation results for models with  regression parameters which  differ from one phase to the other:  $\ef^0_1=(1,0,4,0,-3,5,$ $6,0,-1,0)$, $\ef^0_2=(0,3,-4,-3,0,1,2,-3,0,10)$, $\ef^0_3=(1,3,4,0,0,1,0,0,0,1)$. The random vector $\XX$ is as in the Subsection 5.1.\\
\hh  In Tables \ref{Tabl6} and \ref{Tabl6bis}, we present simulation results when  the regression parameters for the first two phases are the same ($\ef^0_1 = \ef^0_2$) and only the error  distributions are different. In Tables \ref{Tabl5} and  \ref{Tabl6} we present the percentage of true  and of false zero, in each interval, for the three methods, when different error distributions are in each phase. In Tables \ref{Tabl5bis} and  \ref{Tabl6bis} we give the median of the change-points estimations by the three methods. In order to simply the presentation, in these four tables, it was noted by \textit{aQ} the adapted LASSO quantile method, by \textit{Lt} the LASSO-type method for median model and by \textit{aLS} the adaptive LASSO method for LS model.  We calculate bias of the each estimates, the mean of the differences $(\hat \ef-\efo)_{{\cal A}^0}$ and the approximation of the estimation variances $1/M \|(\hat \ef-\efo)_{{\cal A}^0}\|^2_2$, with $M$ the Monte Carlo replications number, and with ${\cal A}^0$ the index set of the nonzero true values. In Tables \ref{Tabl5} - \ref{Tabl6bis}, the law ${\cal E}_1$ is  ${\cal E}xp(-4.5,1)$, ${\cal E}_2$ is  ${\cal E}xp(1.5,1)$,  ${\cal E}_3$ is ${\cal E}xp(-6.5,1)$. The  Cauchy distribution is  ${\cal C}(0,1)$ and Gaussian distribution is ${\cal N}(0,1)$. \\

\noindent \textit{Comparison of the obtained results by the three estimation methods.}\\
\hh  For phases where error distributions are exponential  or Cauchy, the adaptive LASSO quantile method gives the best results in terms  of true zeros  or false zeros percentage, bias and precision of the nonzero parameters.\\
\hh Note that, in all situations when the regression parameters are different from one phase to the other, by the three methods, the median of the change-points estimations coincides or is very close to the true value.\\
\hh  Comparing the last three rows of the Tables \ref{Tabl5bis} and  \ref{Tabl6bis}, we conclude that the adaptive LASSO method for LS model gives less accurate change-points estimates when the change-point is between two phases with the same regression parameters and  error distributions are two exponential law.    Generally, for the LS model with adaptive LASSO penalty, satisfactory results are  obtained in a phase with Gaussian errors, while for Exponential or Cauchy distributions,  poor results are obtained. Recall that for a model with a single phase, the adaptive LASSO method for LS model gave satisfactory results for Gaussian and Exponential errors, and poor results for Cauchy errors (see Tables \ref{Tabl1}, \ref{Tabl2}, \ref{Tabl3}). \\
\hh  By  LASSO-type method, the phases with Exponential errors are poorly estimated.\\
\hh  Comparing the Tables \ref{Tabl5} with \ref{Tabl6} and Tables \ref{Tabl5bis} with \ref{Tabl6bis}, we deduce that selection percentage of the true zeros by adaptive LASSO quantile method decreases slightly for the first phase when $\ef^0_1= \ef^0_2$. However,   this method is better than  two other methods, especially with regards to the bias and precision of estimators  for the nonzero regression parameters. \\

\noindent {\it Conclusion}\\
\hh To conclude, by the adaptive LASSO quantile and LASSO-type methods, the change-points estimation did not affect the sparsity property of the regression parameter estimates. On the other hand, the adaptive LASSO quantile estimators for the nonzero regression parameters are more accurate (in terms of bias and variance) than corresponding LASSO-type estimators.

\begin{table}[p]
{\scriptsize
\caption{\footnotesize Model with three phases.  $\ef^0_1 \neq \ef^0_2 \neq \ef^0_3$. Percentage true 0 and of false 0 by adaptive LASSO quantile, LASSO-type for median model and adaptive LASSO for LS model.}
\begin{center}
\rotatebox{90}{
\begin{tabular}{|c|c|c|c|c|c|c|c|c|c|c|c|c|c|c|c|c|c|c|}\hline
    error distribution & \multicolumn{6}{c|}{interval $(1, \hat l_1)$} &  \multicolumn{6}{c|}{interval $(\hat l_1, \hat l_2)$} &  \multicolumn{6}{c|}{interval $(\hat l_2, n)$}   \\ 
    \cline{2-19} 
     & \multicolumn{3}{c|}{$\%$ true 0} &  \multicolumn{3}{c|}{$\%$ false 0} & \multicolumn{3}{c|}{$\%$ true 0} &  \multicolumn{3}{c|}{$\%$ false 0} & \multicolumn{3}{c|}{$\%$ true 0} &  \multicolumn{3}{c|}{$\%$ false 0} \\
     &  aQ & Lt & aLS &  aQ & Lt & aLS &  aQ & Lt & aLS &  aQ & Lt & aLS &  aQ & Lt & aLS &  aQ & Lt & aLS \\ \hline
     $\varepsilon_1, \varepsilon_2, \varepsilon_3 \sim {\cal E}_1$ & 0.96 & 0.85 & 0.82 & 0.01 & 0.11 & 0.11 & 0.99 & 0.78 & 0.92 & 0 & 0.003 & 0.02 & 1 & 0.67 & 1& 0 & 0.004 & 0.02 \\ \hline 
$\varepsilon_1, \varepsilon_2, \varepsilon_3 \sim {\cal N}$ & 0.97 & 0.97 & 0.88 & 0.01 & 0.04 & 0.12 & 0.98 & 0.98 & 0.96 & 0 & 0 & 0.03 & 0.99 & 0.98 & 1& 0 & 0  & 0.008 \\ \hline 
 $\varepsilon_1, \varepsilon_2, \sim {\cal E}_1, \varepsilon_3 \sim {\cal N}$ & 0.97 & 0.84 & 0.79 & 0.01 & 0.10 & 0.11 & 0.99 & 0.77 & 0.92 & 0 & 0.003 & 0.02 & 0.998 & 0.99 & 1& 0 & 0  & 0.02 \\ \hline 
$\varepsilon_1, \varepsilon_2, \sim {\cal N}, \varepsilon_3 \sim  {\cal E}_1$ & 0.97 & 0.97 & 0.85 & 0.02 & 0.04 & 0.11 & 0.99 & 0.99 & 0.96 & 0 & 0 & 0.03 & 1 & 0.66 & 1& 0 & 0 .004 & 0.03 \\ \hline 
$\varepsilon_1, \varepsilon_3, \sim {\cal E}_1, \varepsilon_2 \sim {\cal N}$ & 0.97 & 0.86 & 0.83 & 0.02 & 0.11 & 0.1 & 0.99 & 0.99 & 0.97 & 0 & 0 & 0.02 & 1 & 0.65 & 1 & 0 & 0.008  & 0.02 \\ \hline 
  $\varepsilon_1, \varepsilon_3, \sim {\cal N}, \varepsilon_2 \sim {\cal E}_1$ & 0.96 & 0.97 & 0.87 & 0.02 & 0.04 & 0.08 & 0.995 & 0.80 & 0.90 & 0  & 0.005 & 0.02 & 0.997 & 0.99 & 1 & 0 & 0  & 0.04 \\ \hline 
  $\varepsilon_1, \varepsilon_3, \sim {\cal E}_1, \varepsilon_2 \sim {\cal C}$ & 0.97 & 0.845 & 0.72 & 0.01 & 0.12 & 0.12 & 0.81 & 0.88 & 0.22 & 0.03 & 0.04 & 0.06 & 0.998 & 0.64 & 0.96 & 0 & 0.007  & 0.05 \\ \hline 
$\varepsilon_1, \varepsilon_3, \sim {\cal E}_2, \varepsilon_2 \sim {\cal E}_1$ & 0.97 & 0.91 & 0.81 & 0.01 & 0.09 & 0.10 & 0.99 & 0.80 & 0.92 & 0 & 0.005 & 0.002 & 0.998 & 0.69 & 1 & 0 & 0  & 0.03 \\ \hline 
  $\varepsilon_1, \varepsilon_3, \sim {\cal E}_3, \varepsilon_2 \sim {\cal E}_1$ & 0.97 & 0.78 & 0.75 & 0.02 & 0.16 & 0.09 & 1 & 0.76 & 0.92 & 0 & 0.005 & 0.02 & 1 & 0.55 & 0.99 & 0 & 0.03  & 0.05 \\ \hline 
   $\varepsilon_1 \sim {\cal E}_3, \varepsilon_2,  \varepsilon_3 \sim {\cal E}_1$ & 0.97 & 0.77 & 0.74 & 0.08 & 0.15 & 0.1 & 0.997 & 0.81 & 0.93 & 0 & 0.007 & 0.02 & 0.998 & 0.65 & 1 & 0 & 0.004  & 0.04 \\ \hline 
\end{tabular}
}
 \end{center}
\label{Tabl5} 
}
\end{table}

%\begin{table}[p]
\begin{table} 
{\scriptsize
\caption{\footnotesize Model with three phases.  $\ef^0_1 \neq \ef^0_2 \neq \ef^0_3$. Summary statistics by adaptive LASSO quantile, LASSO-type for median model and adaptive LASSO for LS model. }
\begin{center}
%\rotatebox{90}{
\begin{tabular}{|c|c|c|c|c|c|c|c|c|c|c|c|c|c|c|c|}\hline
   error distribution      & \multicolumn{3}{c|}{median($\hat l_1$)} &  \multicolumn{3}{c|}{median($\hat l_2$)} & \multicolumn{3}{c|}{mean$(\hat \ef-\efo)_{{\cal A}^0}$} &  \multicolumn{3}{c|}{mean$|(\hat \ef-\efo)_{{\cal A}^0}|$} & \multicolumn{3}{c|}{$1/M\|(\hat \ef-\efo)_{{\cal A}^0}\|^2_{2}$} \\
     &  aQ & Lt & aLS &  aQ & Lt & aLS &  aQ & Lt & aLS &  aQ & Lt & aLS &  aQ & Lt & aLS  \\ \hline
     $\varepsilon_1, \varepsilon_2, \varepsilon_3 \sim {\cal E}_1$ & 30 & 30 & 30 & 100 & 100 & 100 & -0.03 & -0.25 & -0.31 & 0.15 &  0.49 & 0.65 & 0.90 & 7.5 & 10.3 \\ \hline 
$\varepsilon_1, \varepsilon_2, \varepsilon_3 \sim {\cal N}$ & 30 & 31 & 31 & 100 & 100 & 100 & -0.02 & -0.03 & -0.31 & 0.17 &  0.17 & 0.65 & 1 & 1.12 & 9.8 \\ \hline 
 $\varepsilon_1, \varepsilon_2, \sim {\cal E}_1, \varepsilon_3 \sim {\cal N}$ &30 & 30 & 30  &100 & 100 & 100 & -0.03 & -0.17 & -0.31 & 0.16 & 0.26 & 0.64 & 0.94 & 2.87 & 9.9\\ \hline 
  $\varepsilon_1, \varepsilon_3, \sim {\cal E}_1, \varepsilon_2 \sim {\cal N}$ & 31 & 30 & 31 & 100 & 100 & 100 & -0.02 & -0.22 & -0.32 & 0.16 &  0.38 & 0.64 & 0.93 & 5.7 & 9.8 \\ \hline 
  $\varepsilon_1, \varepsilon_3, \sim {\cal N}, \varepsilon_2 \sim {\cal E}_1$ & 30 & 30 & 30 & 100 & 100 & 100 & -0.03 & -0.08 & -0.29 & 0.16 &  0.28 & 0.64 & 0.96 & 3 & 9.5 \\ \hline 
   $\varepsilon_1, \varepsilon_3, \sim {\cal E}_1, \varepsilon_2 \sim {\cal C}$ & 31 & 31 & 31 & 100 & 100 & 100 & -0.05 & -0.25 & -0.37 & 0.26 &  0.49 & 4.8 & 2.7 & 7.5 & 41478 \\ \hline 
  $\varepsilon_1, \varepsilon_3, \sim {\cal E}_2, \varepsilon_2 \sim {\cal E}_1$ & 30 & 30 & 30 & 100 & 100 & 100 & -0.02 & 0.03 & -0.28 & 0.15 &  0.40 & 0.64 & 0.84 & 5.08 & 10.01 \\ \hline
 $\varepsilon_1, \varepsilon_3, \sim {\cal E}_3, \varepsilon_2 \sim {\cal E}_1$ & 30 & 30 & 31 & 100 & 100 & 100 & -0.03 & -0.34 & -0.34 & 0.15 &  0.62 & 0.68 & 0.87 & 13.2 & 11.9 \\ \hline 
  $\varepsilon_1 \sim {\cal E}_3, \varepsilon_2,  \varepsilon_3 \sim {\cal E}_1$ & 30 & 30 & 31 & 100 & 100 & 100 & -0.02 & -0.31 & -0.34 & 0.15 &  0.56 & 0.68 & 0.81 & 10.7 & 11.6 \\ \hline
\end{tabular}
%}
 \end{center}
\label{Tabl5bis} 
}
\end{table}

\begin{table}[p]
{\scriptsize
\caption{\footnotesize Model with three phases.  $\ef^0_2 = \ef^0_1$,  $\ef^0_2\neq \ef^0_3$. Percentage true 0 and of false 0 by adaptive LASSO quantile, LASSO-type for median model and adaptive LASSO for LS model.}
\begin{center}
\rotatebox{90}{
\begin{tabular}{|c|c|c|c|c|c|c|c|c|c|c|c|c|c|c|c|c|c|c|}\hline
    error distribution & \multicolumn{6}{c|}{interval $(1, \hat l_1)$} &  \multicolumn{6}{c|}{interval $(\hat l_1, \hat l_2)$} &  \multicolumn{6}{c|}{interval $(\hat l_2, n)$}   \\ 
    \cline{2-19} 
     & \multicolumn{3}{c|}{$\%$ true 0} &  \multicolumn{3}{c|}{$\%$ false 0} & \multicolumn{3}{c|}{$\%$ true 0} &  \multicolumn{3}{c|}{$\%$ false 0} & \multicolumn{3}{c|}{$\%$ true 0} &  \multicolumn{3}{c|}{$\%$ false 0} \\
     &  aQ & Lt & aLS &  aQ & Lt & aLS &  aQ & Lt & aLS &  aQ & Lt & aLS &  aQ & Lt & aLS &  aQ & Lt & aLS \\ \hline
 $\varepsilon_1, \varepsilon_3, \sim {\cal E}_1, \varepsilon_2 \sim {\cal N}$ & 0.87 & 0.86 & 0.84 & 0.07 & 0.11 & 0.12 & 0.98 & 0.99 & 0.98 & 0.001 & 0 & 0.05 & 1 & 0.67 & 0.996& 0 & 0.002  & 0.04 \\ \hline 
  $\varepsilon_1, \varepsilon_3, \sim {\cal N}, \varepsilon_2 \sim {\cal E}_1$ & 0.90 & 0.95 & 0.84 & 0.05 & 0.03 & 0.11 & 0.98 & 0.78 & 0.97 & 0.001 & 0.01 & 0.05 & 0.99 & 0.98 & 1 & 0 & 0  & 0.16 \\ \hline 
  $\varepsilon_1, \varepsilon_3, \sim {\cal E}_1, \varepsilon_2 \sim {\cal C}$ & 0.88 & 0.91 & 0.74 & 0.08 & 0.12 & 0.11 & 0.91 & 0.94 & 0.49 & 0.01 & 0.01 & 0.08 & 1 & 0.64 & 0.97 & 0 & 0  & 0.05 \\ \hline 
  $\varepsilon_1, \varepsilon_3, \sim {\cal E}_2, \varepsilon_2 \sim {\cal E}_1$ & 0.81 & 0.88 & 0.83 & 0.1 & 0.08 & 0.12 & 0.99 & 0.80 & 0.95 & 0.002 & 0.02 & 0.06 & 1 & 0.69 & 0.99 & 0 & 0  & 0.02 \\ \hline 
  $\varepsilon_1, \varepsilon_3, \sim {\cal E}_3, \varepsilon_2 \sim {\cal E}_1$ & 0.94 & 0.78 & 0.82 & 0.03 & 0.19 & 0.12 & 0.99 & 0.78 & 0.96 & 0 & 0.01 & 0.05 & 0.99 & 0.56 & 0.99 & 0 & 0.04  & 0.05 \\ \hline 
   $\varepsilon_1 \sim {\cal E}_3, \varepsilon_2,  \varepsilon_3 \sim {\cal E}_1$ & 0.92 & 0.77 & 0.78 & 0.04 & 0.18 & 0.11 & 1 & 0.8 & 0.95 & 0 & 0.01 & 0.04 & 1 & 0.64 & 1 & 0 & 0  & 0.03 \\ \hline 
\end{tabular}
}
 \end{center}
\label{Tabl6} 
}
\end{table}

\begin{table} 
\caption{\footnotesize Model with three phases.  $\ef^0_2 =\ef^0_1$, $\ef^0_2 \neq \ef^0_3$. Summary statistics by adaptive LASSO quantile, LASSO-type for median model and adaptive LASSO for LS model, $n=200$, $l^0_1=30$, $l^0_2=100$. }
\begin{center}
{\scriptsize
%\rotatebox{90}{
\begin{tabular}{|c|c|c|c|c|c|c|c|c|c|c|c|c|c|c|c|}\hline
   error distribution      & \multicolumn{3}{c|}{median($\hat l_1$)} &  \multicolumn{3}{c|}{median($\hat l_2$)} & \multicolumn{3}{c|}{mean$(\hat \ef-\efo)_{{\cal A}^0}$} &  \multicolumn{3}{c|}{mean$|(\hat \ef-\efo)_{{\cal A}^0}|$} & \multicolumn{3}{c|}{$1/M\|(\hat \ef-\efo)_{{\cal A}^0}\|^2_{2}$} \\
     &  aQ & Lt & aLS &  aQ & Lt & aLS &  aQ & Lt & aLS &  aQ & Lt & aLS &  aQ & Lt & aLS  \\ \hline
   $\varepsilon_1, \varepsilon_3, \sim {\cal E}_1, \varepsilon_2 \sim {\cal N}$ & 32 & 31 & 28 & 100 & 100 & 100 & -0.04 & -0.23 & -0.33 & 0.15 &  0.43 & 0.63 & 0.79 & 6.1 & 8.8 \\ \hline 
  $\varepsilon_1, \varepsilon_3, \sim {\cal N}, \varepsilon_2 \sim {\cal E}_1$ & 28 & 30 & 26 & 100 & 100 & 100 & -0.04 & -0.12 & -0.32 & 0.22 &  0.29 & 0.65 & 2.3 & 3.06 & 9.6 \\ \hline 
   $\varepsilon_1, \varepsilon_3, \sim {\cal E}_1, \varepsilon_2 \sim {\cal C}$ & 31 & 31 & 39 & 100 & 100 & 99 & -0.06 & -0.22 & -0.35 & 0.28 &  0.44 & 1.32 & 3.2 & 6.1 & 316 \\ \hline 
  $\varepsilon_1, \varepsilon_3, \sim {\cal E}_2, \varepsilon_2 \sim {\cal E}_1$ & 32 & 31 & 30 & 100 & 100 & 100 & -0.07 & -0.007 & -0.29 & 0.25 &  0.42 & 0.64 & 3.2 & 5.2 & 9.2 \\ \hline
 $\varepsilon_1, \varepsilon_3, \sim {\cal E}_3, \varepsilon_2 \sim {\cal E}_1$ & 31 & 30 & 25 & 101 & 100 & 100 & -0.04 & -0.45 & -0.38 & 0.19 &  0.71 & 0.72 & 1.66 & 16.9 & 12.45 \\ \hline 
  $\varepsilon_1 \sim {\cal E}_3, \varepsilon_2,  \varepsilon_3 \sim {\cal E}_1$ & 34 & 30 & 27 & 100 & 100 & 100 & -0.05 & -0.37 & -0.38 & 0.22 &  0.63 & 0.74 & 2.9 & 13.4 & 13 \\ \hline
\end{tabular}
}
\end{center}
\label{Tabl6bis} 
\end{table}

\subsubsection{Estimation of the change-point number}
In  Table \ref{Tabl7} we give results, after 100  Monte Carlo replications, in order to estimate  the  number of  phases using relation (\ref{eq19}). There was one change-point to the observation $l^0_1=30$ for a  total of 100 observations. The change in the model is due either to the change in the regression parameters ($\ef_1^0 \neq \ef_2^0$) or to the change in quantile of the error (for the same  index $\tau$).\\
\hh We compare the criterion  proposed  in this paper for the  adaptive LASSO quantile  method with the criterion proposed in the paper of  \cite{Ciuperca:13a}, for the  adaptive LASSO method for LS model. In the paper of \cite{Ciuperca:13b}, where the  LASSO-type method for a median  model with change-points has been studied, there is no criterion proposed to estimate the true number of change-points. Therefore, we propose here a criterion of the same shape as for the  adaptive LASSO quantile method: to the corresponding penalized objective function add a term of the form $G(K,p)B_n$ with $G(K,p)=K$  and the sequence $B_n=n^{5/8}$. The criterion values for the three methods are calculated for $K \in \{0,1,2,3 \}$.    For the adaptive LASSO quantile method we consider $\tau=0.55$. \\
\hh From Table \ref{Tabl7}, we deduce that,  the criteria associated to the three methods choose correctly  the change-point number when the change is due to the regression  parameters. On the other hand, when the change is due only to the quantile, i.e. $\ef^0_1=\ef^0_2=(1,0,4,0,-3,5,6,0,-1,0)$, then, if the errors come from the same distribution (Exponential) but with different quantiles, the criterion for the adaptive LASSO quantile method does not identify whenever the change (only 62/100), the criterion for the  LASSO-type method never identifies the change, while that for the adaptive LASSO method  for LS model identifies 10/100. The results improve  for adaptive LASSO quantile and LASSO-type methods when the errors of the two phases have different distributions  (Exponential and Gaussian). 

\begin{table}
{\scriptsize
\caption{\footnotesize Results on the choice of the change-point number by the criteria associated to the  methods: adaptive LASSO quantile, adaptive LASSO for LS model and LASSO-type for median(LAD) model, $K^0=1$. 100 Monte Carlo replications. Number of $\hat K$ = 0,1,2,3 for 100 Monte Carlo replications. }
\begin{center}
\begin{tabular}{|c|c|c|c|c|c|c|c|c|c|c|c|c|c|}\hline
$\ef_1, \ef_2$ &   error distribution      & \multicolumn{4}{c|}{QUANT+aLASSO} &  \multicolumn{4}{c|}{LS+aLASSO} & \multicolumn{4}{c|}{LAD+LASSOtype}  \\ 
 & & \multicolumn{4}{c|}{number of $\hat K=$} & \multicolumn{4}{c|}{number of $\hat K=$} & \multicolumn{4}{c|}{number of $\hat K=$} \\
  & & 0 & 1 & 2 & 3 & 0 & 1 & 2 & 3 & 0 & 1 & 2 & 3 \\ 
  \hline
$\ef^0_1=\ef^0_2$  & $\varepsilon_1 \sim {\cal E}_1, \varepsilon_2 \sim {\cal E}_3$  & 36 & 62 & 1 & 1& 88 & 10 & 1 & 1 & 100 & 0 & 0 & 0 \\  
 & $\varepsilon_1 \sim {\cal E}_1, \varepsilon_2 \sim {\cal N}(0,1)$  & 0 & 99 & 1 & 0 & 54 & 66 & 0 & 0 & 1 & 99 & 0 & 0 \\  \hline 
  $\ef^0_1 \neq \ef^0_2$  & $\varepsilon_1 \sim {\cal E}_1, \varepsilon_2 \sim {\cal E}_3$  & 0 & 100 & 0 & 0& 0 & 100 & 0 &0 & 0 & 100 & 0 & 0 \\  
  & $\varepsilon_1 \sim {\cal E}_1, \varepsilon_2 \sim {\cal N}(0,1)$  &  0 & 100 & 0 & 0& 0 & 100 & 0 &0 & 0 & 100 & 0 & 0 \\ \hline   
\end{tabular}
 \end{center}
\label{Tabl7} 
}
\end{table}
%\clearpage
 \newpage
\section{Proofs} 
For the convenience of the reader, recall first a lemma due to \cite{Babu:89}.
\begin{lemma}
\label{l 1} (\cite{Babu:89}, Lemma 1)\\
Let $Z_i$ be a sequence of independent random variables with mean zero and $|Z_i| \leq \beta $ for some $\beta >0$. Let also $V \geq \sum^n_{i=1}\eE[Z_i^2] $. Then for all $0<s<1$ and $0 \leq z \leq V/(s \beta)$, we have
\begin{equation}
\label{e02}
\eP\cro{\left| \sum^n_{i=1}Z_i \right| >z} \leq 2 \exp \pth{-z^2s(1-s)/V}.
\end{equation}
\end{lemma}

This section is divided into two subsections. In the first subsection we give the proofs of all Theorems and Propositions. In the second subsection, we present the lemma proofs.

\subsection{Proposition and Theorem proofs}
\noindent {\bf Proof of Proposition \ref{KKT}}
\textit{(i)} For all $j \in {\cal A}^*_n$, the estimator  $\hat \phi^*_{n,j}$ is the solution of the equation:
$
0=\sum^n_{i=1} \frac{\partial R^{(\tau,\lambda)}_i (b, \ef; b^0, \efo)}{\partial \phi_{,j}}
$. 
By elementary algebra  we obtain that 
\[ \sum^n_{i=1} \frac{\partial R^{(\tau,\lambda)}_i (b, \ef; b^0, \efo)}{\partial \phi_{,j}} = - \tau X_{ij} +X_{ij} \e1_{Y_i < \XX^t_i \ef}+ \lambda_n \hat \omega_{n,j}sgn(\phi_{,j})\]  and the assertion affirmation \textit{(i)} follows.\\
\textit{(ii)} In this case, the subgradient set $\partial \| 0_j \|_1$ is the closed interval $[-1,1]$. Then 
\[
0 \in \sum^n_{i=1} \frac{\partial R^{(\tau,\lambda)}_i (b, \ef; b^0, \efo)}{\partial \phi_{,j}} =- \tau X_{ij} +X_{ij} \e1_{Y_i < \XX^t_i \ef}+ \lambda_n \hat \omega_{n,j} [-1,1].\]
This leads  conclusion.
 \hspace*{\fill}$\blacksquare$ \\
 
 \noindent {\bf Proof of Theorem \ref{Theorem 4.1 ZouYuan}} \textit{(i)} We reparameterize the model: $\uu_n =\sqrt{n} (\hat \ef^*_n - \efo)$ and $u_{0,n} = \sqrt{n} (\hat b^*_n -b^0)$. Then $(u_{0,n}, \uu_n)$ is the minimizer of the criterion 
\[
\sum^n_{i=1} \cro{\rho_\tau \pth{\varepsilon_i-b^0-n^{-1/2}(u_{0,n}+\XX^t_i \uu_n)} - \rho_\tau (\varepsilon_i - b^0)}  +\lambda_n \hat \eo^t_n [ | \efo +\uu_n n^{-1/2}  | -|\efo|].\]
The rest of the proof is similar as that of Theorem 4.1 in the paper of  \cite{Zou:Yuan:08} and we omit it. With the remark that in our case (more general than in \cite{Zou:Yuan:08}) the supposition $n^{(g-1)/2} \lambda_n \rightarrow \infty$ is required to prove, in the case $\phi^0_{,j}=0$ and $u_{n,j} \neq 0$, that 
\[
\frac{\lambda_n}{\sqrt{n} |\hat \ef_n|^g} \sqrt{n} \pth{\left| \phi^0_{,j}+\frac{u_{n,j}}{\sqrt{n}}  \right| -|\phi^0_{,j}| } \overset{\eP} {\underset{n \rightarrow \infty}{\longrightarrow}} \infty .\]
\textit{(ii)}  If $j \in {{\cal A}^0}$, then $\phi^0_{,j} \neq 0$. Since $\hat \phi^*_{n,j}$ is the corresponding estimator of $\phi^0_{,j}$, and is asymptotically normal, thus $j  \in \hat {\cal A}^*_n$ with probability  tending to 1. Then $\lim_{n \rightarrow \infty} \eP[\hat {\cal A}^*_n \supseteq {{\cal A}^0} ]=1$. \\
For the reverse inclusion, we show that if $j \not \in {{\cal A}^0}$ then $j \not \in \hat {\cal A}^*_n$. Let us calculate $\eP[j  \in  \hat  {\cal A}^*_n, j \not \in {\cal A}^0]$. Since $j  \in {\cal A}^*_n$ it follows from the KKT optimality conditions  of the Proposition \ref{KKT}(i), that
\[\tau \sum^n_{i=1} X_{ij} -  \sum^n_{i=1} X_{ij}  \e1_{Y_i < \XX^t_i \hat \ef^*_n}= \lambda_n \hat \omega_{n,j} sgn(\hat \phi^*_{n,j}),\]
 which implies that  
\begin{equation}
\label{eqZ}
\lambda_n \hat \omega_{n,j} < 2 \sum^n_{i=1} | X_{ij}|.
\end{equation}
For the left member of the last inequality, we have
\begin{equation}
\label{eqB}
\frac{\lambda_n \hat \omega_{n,j}}{n}=\lambda_n \frac{1}{|n^{1/2} \hat \phi_{n,j} |^g} \cdot \frac{n^{g/2}}{n},
\end{equation}
where $ \hat \phi_{n,j}$ is the quantile estimator of $\phi^0_{,j}$. On the other hand, since $j \not \in {\cal A}^0 $, we have that $\phi^0_{,j} =0$. Using the fact that $ \hat \phi_{n,j}$ is strongly consistent and asymptotically normal, we have that for all $\epsilon>0$, there exists $\eta_\epsilon >0$ such that 
\begin{equation}
\label{eqC}
\eP[n^{-1/2} | \hat \phi_{n,j}|^{-1} > \eta_\epsilon ] >1 -\epsilon .
\end{equation} 
Since  $n^{g/2-1} \lambda_n \rightarrow \infty$, with the relation (\ref{eqC}), we obtain that (\ref{eqB}) converges to infinity with probability converging to 1 as $n \rightarrow \infty$. On the other hand, by the Cauchy-Schwarz inequality
$n^{-1} \sum^n_{i=1} | X_{ij}| \leq \pth{n^{-1} \sum^n_{i=1} X^2_{ij} }^{1/2} $, which is bounded with probability converging 1, by assumption (A1). Taking into account (\ref{eqZ}), we obtain that  (\ref{eqB}) is bounded. Contradiction. This completes the proof.
 \hspace*{\fill}$\blacksquare$ \\
 
 \noindent {\bf Proof of Proposition \ref{proposition 2.1.SCAD}}
Let us denote $h_i=b^0-b+\XX^t_i(\efo-\ef)$. Then 
\[\eE[R^{(\tau)}_i (b,\ef;b^0,\efo)]=\int_{\R} [\rho_\tau(u+h_i)-\rho_\tau(u)] dF(u+b^0).\]
\hh \underline{If $h_i \leq 0 $}. By elementary calculations we obtain
$
\int^0_{- \infty} [\rho_\tau(x-h_i)-\rho_\tau(x)] dF(x+b^0)=-\tau(1-\tau)h_i
$
and also
$
\int^{\infty}_0 [\rho_\tau(x-h_i)-\rho_\tau(x)] dF(x+b^0)= \tau(1-\tau)h_i- \int^{-h_i}_0 [h_i(1-\tau)+\tau h_i +x] dF(x+b^0) $. Then, in this case
$
\eE[R^{(\tau)}_i (b,\ef;b^0,\efo)] =\int^{-h_i}_0[-x+ |h_i|] dF(x+b^0)$.\\
\hh \underline{If $h_i > 0 $}. By similar calculations as above we obtain that 
$
\eE[R^{(\tau)}_i (b,\ef;b^0,\efo)]=\int_{-h_i}^0[x+ |h_i|] dF(x+b^0)$.\\
Thus, we can write all this in a more condensed rule:
\[
\eE[R^{(\tau)}_i (b,\ef;b^0,\efo)] \geq \e1_{h_i >0} \frac{|h_i|}{2}\int_{-\frac{h_i}{2}}^0 dF(x+b^0)+ \e1_{h_i \leq 0} \frac{|h_i|}{2}\int^{-\frac{h_i}{2}}_0 dF(x+b^0)
\]
\[
\qquad \qquad \qquad \qquad = \frac{|h_i|}{2} \cro{\e1_{h_i >0}  [F(b^0)-F(b^0-h_i/2)] +  \e1_{h_i \leq 0} [F(b^0-h_i/2)- F(b^0)]} \geq 0.
\]
 \hspace*{\fill}$\blacksquare$ \\

    \noindent {\bf Proof of Theorem \ref{Theorem 5aLASSO}. }
  \textbf{\underline{Step I.}} We prove that, with probability approaching 1, the adaptive LASSO quantile change-point estimators are to a smaller distance than $[n^{1/2}]$. For this purpose,  let us consider the constant $\varrho \in (\alpha, \min(1, (1+g)/2))$, with $\alpha>1/2$ as in Lemma \ref{lemma 6aLASSO}. We will prove that: $\eP[|\hat l^*_r-l^0_r | >n^\varrho] \rightarrow 0$ as $n \rightarrow \infty$, for each $r=1, \cdots, K$. \\
  Consider the set of change-points, all close to the true points at a distance less $[n^\varrho]$:
  \[{\cal L}(\varrho) \equiv \{(l_1, \cdots, l_K); \sum^K_{r=} |l_r-l^0_r| \leq [n^\varrho] \}.\]
   Consider a subset of its complementary, for some $r \in \{ 1, \cdots, K\}$:
   \[{\cal L}^c_r(\varrho) \equiv \{(l_1, \cdots, l_K); |l_t-l^0_r| >n^\varrho, \forall t = 1, \cdots, K   \},\] with all the change-points to a distance of $l^0_r$ greater than $[n^\varrho]$. By the definition of the objective function $S^*$ we have, for all $(l_1, \cdots , l_K) \in {\cal L}^c_r(\varrho)$, with probability 1:
 \begin{equation}
 \label{ineq_10aLASSO}
 S^*(l_1, \cdots, l_K) \geq S^*(l_1,\cdots, l_K, l^0_1, \cdots , l^0_{r-1}, l^0_r-[n^\rho], l^0_r+[n^\rho], l^0_{r+1}, \cdots, l^0_K )  \equiv \sum^{K+2}_{q=1}{\cal T}_q,
 \end{equation}
where  ${\cal T}_q$  are, for $q \in \{1, \cdots, r-1, r+1, \cdots, K+2 \}$, the penalized sums involving observations between $l^0_{q-1}$ and $l^0_q$. ${\cal T}_{r}$  
  is the penalized sum involving observations between $l^0_{r-1}$ and $l^0_r -[n^\varrho]$ and ${\cal T}_{r+1} $ between $l^0_r+[n^\varrho]$ and $l^0_{r+1}$; ${\cal T}_{K+2}$ is calculated between $l^0_r -[n^\varrho]$ and $l^0_r+[n^\varrho]$. Note that the sum $S^*(l_1,\cdots, l_K, l^0_1, \cdots , l^0_{r-1}, l^0_r-[n^\rho], l^0_r+[n^\rho], l^0_{r+1}, \cdots, l^0_K )$ is the extension of the definition (\ref{def_S*}) for $2K+1$ change-points: $l_1,\cdots, l_K, l^0_1, \cdots , l^0_{r-1},$ $ l^0_r-[n^\rho], l^0_r+[n^\rho], l^0_{r+1}, \cdots, l^0_K$.\\ 
  Since $(l_1, \cdots , l_K) \in {\cal L}^c_r(\varrho)$, all change-points $l_1, \cdots, l_K$ are in ${\cal T}_1, \cdots, {\cal T}_{r-1}, \cdots, {\cal T}_{K+1}$ and none  in ${\cal T}_{K+2}$.  For  each $q \in \{1, \cdots, r-1, r+1, \cdots, K+2 \}$, let us consider the points  between two true consecutive change-points $k_{1,q} < \cdots < k_{J(q),q}\equiv \{l_1, \cdots, l_K \} \cap \{ j; l^0_{q-1} < j \leq l^0_q \}$. The number $ J(q)$ is greater (or equal) than zero and smaller (or equal) than $K$, with the property that  for $q \neq q'$, $J(q) \neq J(q')$. The penalized sum ${\cal T}_q$ can be written
\[
  {\cal T}_q=\sum^{J(q)+1}_{j=1} \min_{b_j, \ef_j} [ \sum^{k_{j,q}}_{i=k_{j-1},q+1} \rho_\tau(\varepsilon_i-b_j-\XX^t_i(\ef_j -\ef^0_q))+\lambda_{(k_{j-1,q};k_{j,q})} \hat \eo^t_{(k_{j-1,q};k_{j,q})} | \ef_j|]
.\]
 It is obvious that
 \[
0 \geq {\cal T}_q -  \sum^{J(q)+1}_{j=1} [\sum^{k_{j,q}}_{i=k_{j-1},q+1} \rho_\tau(\varepsilon_i-b^0_q)+\lambda_{(k_{j-1,q};k_{j,q})} \hat \eo^t_{(k_{j-1,q};k_{j,q})} | \ef^0_t| ]\]
\[
\geq  -2(K+1) \sup_{0 \leq l <k \leq n,k-l=d_n} | \inf_{b,\ef} \sum^k_{i=l+1} R^{(\tau,\lambda)}_{i,(l;k)} (b,\ef;b^0, \efo) |,\] 
with $l$ and $k$ two observations without any change between they,  $b^0, \efo$ the true values of the parameters on this interval and $d_n$ a convergent sequence to  $\infty$ when  $n$ converges to $\infty$. Due to Lemma \ref{lemma 6aLASSO}, 
 \begin{equation}
 \label{eq13}
2(K+1) \sup_{0 \leq l <k \leq n,k-l=d_n} | \inf_{b,\ef} \sum^k_{i=l+1} R^{(\tau,\lambda)}_{i,(l;k)} (b,\ef;b^0, \efo) | = -\min(O_{\eP}(n^\alpha),o_{\eP}(n^{(1+g)/2})).
 \end{equation}
 For the observations between $l^0_r-[n^\varrho]$ and $ l^0_r+[n^\varrho]$, we have that:
 \[ {\cal T}_r-\sum^{J(r)+1}_{j=1}[ \sum^{k_{j,r}}_{i=k_{j_1,r}+1} \rho_\tau(\varepsilon_i-b^0_r)+\lambda_{(k_{j-1,r};k_{j,r})} \hat \eo^t_{(k_{j-1,r};k_{j,r})} |\ef^0_r |+\lambda_{(k_{j-1,r};k_{j,r})} \hat \eo^t_{(k_{j-1,r};k_{j,r})} |\ef^0_{r+1} | ]\]
 \[  = \sum^{J(r)+1}_{j=1}\min_{b_j,\ef_j} [ \sum^{k_{j,r}}_{i=k_{j_1,r}+1}\rho_\tau(\varepsilon_i-b_j-\XX^t_i(\ef_j -\ef^0_t)) - \rho_\tau(\varepsilon_i-b^0_r) +\lambda_{(k_{j-1,t};k_{j,t})} \hat \eo^t_{(k_{j-1,t};k_{j,t})} | \ef_j|  ] \]  
 \[ \qquad \qquad \qquad \qquad   - \lambda_{(k_{j-1,r};k_{j,r})} \hat \eo^t_{(k_{j-1,r};k_{j,r})} [|\ef^0_r|+ |\ef^0_{r+1} | ] .\]  
 But for $(l_1, \cdots, l_K) \in {\cal L}^c_r(\varrho)$, there is no other change-point between $l^0_r-[n^\varrho]$ and $l^0_r+[n^\varrho]$, so this last relation is in fact
 \begin{equation}
 \label{eq12}
   \min_{b, \ef} [\sum^{l^0_r}_{i=l^0_r-[n^\varrho]+1} R^{(\tau,\lambda)}_{i;(l^0_r-[n^\varrho];l^0_r)}(b,\ef;b^0_r,\ef^0_r) +\sum^{l^0_r+[n^\varrho]}_{i=l^0_r+1}R^{(\tau,\lambda)}_{i;(l^0_r;l^0_r+[n^\varrho])}(b,\ef;b^0_{r+1},\ef^0_{r+1}) ] . 
 \end{equation}
Since to left and to right of $l^0_r$ we have different models, then $(b^0_r, \ef^0_r) \neq (b^0_{r+1},\ef^0_{r+1})$. This implies that for the $(b,\ef)$ which minimizes (\ref{eq12}) we have either  $|b-b^0_r|+\| \ef-\ef^0_r\|_2 >c >0$ or  $|b-b^0_{r+1}|+\| \ef-\ef^0_{r+1}\|_2 >c >0$. Without loss of generality, consider the last case. Using  Lemma \ref{lemma 7aLASSO} for $c_n=c>0$ on the interval $(l^0_r, l^0_r+[n^\varrho])$, we obtain that,  there exists $\epsilon> 0$ such that
 \[c^{-1}[n^\varrho]^{-1} \sum^{l^0_r+[n^\varrho]}_{i=l^0_r+1} R^{(\tau,\lambda)}_{i;(l^0_r;l^0_r+[n^\varrho])}(b,\ef;b^0_{r+1},\ef^0_{r+1}) > \epsilon >0.\]
Thus
 \begin{equation}
 \label{eq14a}
 \inf_{b, \ef}  \sum^{l^0_r+[n^\varrho]}_{i=l^0_r+1} R^{(\tau,\lambda)}_{i;(l^0_r;l^0_r+[n^\varrho])}(b,\ef;b^0_{r+1},\ef^0_{r+1}) \geq O_{\eP}(n^\varrho ) >0 . 
 \end{equation}
Then, taking into account the relations (\ref{ineq_10aLASSO}),   (\ref{eq13}), (\ref{eq12}), (\ref{eq14a}) and the fact that by definition of $S^*$, we have that $S^*(\hat l^*, \cdots, \hat l^*_K) \leq S^*_0$ with probability 1, we obtain that $S^*(l_1, \cdots, l_K)- S^*_0$ is greater than 
 $ \sum^{K+2}_{q=1} {\cal T}_q -S^*_0 $ which is greater than $ -\min(O_{\eP}(n^\alpha),o_{\eP}(n^{(1+g)/2}))+O_{\eP}(n^\varrho) +\sum^{K+2}_{q=1} \sum^{J(q)+1}_{j=1} \lambda_{(k_{j-1,q}, k_{j,q})} \hat \eo_{(k_{j-1,q}, k_{j,q})}$ $
 - \sum^{K+1}_{r=1}\lambda_{(l^0_{r-1}, l^0_r)} \hat \eo^t_{(l^0_{r-1}; l^0_r)} |\ef^0_r|$ $ =O_{\eP}(n^\varrho)+o_{\eP}(n^{(g+1)/2})=O_{\eP}(n^\varrho) >0 $.
Then, for $n \rightarrow \infty$, 
\[ \eP[\min_{(l_1, \cdots, l_K) \in {\cal L} ^c_r(\varrho)} S^*(l_1, \cdots , l_{K}) > S^*_0] \rightarrow 1.\] 
 On the other hand, since $(\hat l^*, \cdots, \hat l^*_K)$ are the change-points estimators, we have that  $ S^*(\hat l^*_1, \cdots, \hat  l^*_K) \leq S^*_0$, with probability 1. These last two relations imply
 \[\eP [ (\hat l^*_1, \cdots, \hat l^*_K) \in {\cal L}^c_r(\varrho)] \rightarrow 0, \textrm{as } n \rightarrow \infty, \textrm{ for every } r=1, \cdots, K.\]
  \textbf{\underline{Step II.}} We show that for all $\nu <1/4$, for every $ r=1, \cdots, K$ we have $\eP[| \hat l^*_r - l^0_r| >n^\nu] \rightarrow 0$, as $n \rightarrow \infty$.\\
  \begin{figure}
  \includegraphics[scale=1]{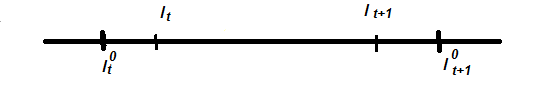}   
  \caption{ }
\label{Figure 1}
  \end{figure}
   By Step I, we have
   $\eP [ (\hat l^*_1, \cdots, \hat l^*_K) \in {\cal L}(\varrho)] \rightarrow 1$, as $n \rightarrow \infty$.  Consider then the change-points $(l_1, \cdots, l_K)$ belonging to the set ${\cal L}(\varrho)$ and $(l_1, \cdots, l_K) \in {\cal L}^c_r(\nu)$, where ${\cal L}^c_r(\nu)$ is similar to the set ${\cal L}^c_r(\varrho)$ with $\nu$ instead of $\varrho$.  For these change-points,  there is a similar relation to (\ref{ineq_10aLASSO}) for the objective function $S^*$. For $q \neq r-1,r$, by  assumption (A3), using  step I, we have that there are at most two points $l_q$ and $l_{q+1}$ between $l^0_q$ and $l^0_{q+1}$. Suppose that there are two points $l_q$ and $l_{q+1}$ between $l^0_q$ and $l^0_{q+1}$ (see Figure \ref{Figure 1}). If there is a single  point or no point the approach is the same.  Let us define the following sums
   \[
 \begin{array}{cl}
 {\cal D}(l^0_q,l^0_{q+1}) \equiv & \inf_{b, \ef} \{  \sum^{l_q}_{i=l^0_q+1} \rho_\tau(\varepsilon_i-b-\XX^t_i(\ef-\ef^0_{q+1})) +\lambda_{(l^0_q; l_q)} \hat \eo^t_{(l^0_q; l_q)}|\ef|  \} \\
 & +\inf_{b, \ef} \{  \sum^{l_{q+1}}_{i=l_q+1} \rho_\tau(\varepsilon_i-b-\XX^t_i(\ef-\ef^0_{q+1})) +\lambda_{(l_q; l_{q+1})} \hat \eo^t_{(l_q; l_{q+1})}|\ef|   \}\\
&  +\inf_{b, \ef} \{  \sum^{l^0_{q+1}}_{i=l_{q+1}+1} \rho_\tau(\varepsilon_i-b-\XX^t_i(\ef-\ef^0_{q+1})) +\lambda_{(l_{q+1}; l^0_{q+1})} \hat \eo^t_{(l_{q+1}; l^0_{q+1})}|\ef|  \} .
\end{array}
\]
The sums $ {\cal D}(l^0_{r-1}, l^0_{r}-[n^\nu])$, $ {\cal D}(l^0_{r}+[n^\nu],l^0_{r+1})$, ${\cal D}( l^0_{r}-[n^\nu],l^0_{r}+[n^\nu])$ can be defined in the same way. 
Then, for the difference between the following two objective functions 
\[S^*(l_1,\cdots, l_K, l^0_1, \cdots , l^0_{r-1}, l^0_r-[n^\nu], l^0_r+[n^\nu], l^0_{r+1}, \cdots, l^0_K )- S^*_0,\] we focus on what happens between two change-points. This last difference can be written
\[
 \begin{array}{l}
\sum_{q\neq r-1,r} \left\{ {\cal D}(l^0_q, l^0_{q+1}) \right.- \left. \sum^{l^0_{q+1}}_{i=l^0_q+1} \rho_\tau(\varepsilon_i-b^0_{q+1}) -\lambda_{(l^0_q; l^0_{q+1})} \hat \eo^t_{(l^0_q; l^0_{q+1})} |\ef^0_{q+1}| \right\}\\
 +\left\{ {\cal D}(l^0_{r-1}, l^0_{r}-[n^\nu])- \sum^{l^0_{r}-[n^\nu]}_{i=l^0_{r-1}+1}\rho_\tau(\varepsilon_i-b^0_{r})\right. \left. - \lambda_{(l^0_{r-1}; l^0_{r}-[n^\nu])} \hat \eo^t_{(l^0_{r-1}; l^0_{r}-[n^\nu])} |\ef^0_{r}| \right\}\\
 +\left\{ {\cal D}(l^0_{r}+[n^\nu],l^0_{r+1})- \sum^{l^0_{r+1}}_{i=l^0_{r}+[n^\nu]+1}\rho_\tau(\varepsilon_i-b^0_{r+1})\right. $  $\left. - \lambda_{(l^0_{r}+[n^\nu]; l^0_{r+1})} \hat \eo^t_{(l^0_{r}+[n^\nu]; l^0_{r+1})} |\ef^0_{r+1}| \right\} \\
 +\left\{ {\cal D}( l^0_{r}-[n^\nu],l^0_{r}+[n^\nu])\right.- \left. \sum^{l^0_{r}+[n^\nu]}_{i=l^0_{r}+1}\rho_\tau(\varepsilon_i-b^0_{r+1}) \right. \left.- \lambda_{(l^0_{r}; l^0_{r}+[n^\nu])} \hat \eo^t_{(l^0_{r}; l^0_{r}+[n^\nu])} |\ef^0_{r+1}| 
\right.  \\
  - \left. \sum^{l^0_{r}}_{i=l^0_{r}-[n^\nu]+1}\rho_\tau(\varepsilon_i-b^0_{r})\right.   \left. - \lambda_{(l^0_{r}-[n^\nu]; l^0_{r})} \hat \eo^t_{(l^0_{r}-[n^\nu]; l^0_{r})} |\ef^0_{r}| \right\} 
 \equiv D_1+D_2+D_3+D_4. \end{array}
\]
By Remark \ref{Remark 1}, we have that $D_1, D_2, D_3=O_{\eP}(1)$. \\
For $D_4$ we have by the definition of sums  $ {\cal D}$ that $ {\cal D}(l^0_{r}-[n^\nu],l^0_{r}+[n^\nu])$ is equal to 
\[
 \begin{array}{l} \qquad
\inf_{b,\ef} \{ \sum^{l^0_{r}}_{i=l^0_{r}-[n^\nu]+1}\rho_\tau(\varepsilon_i-b-\XX^t_i(\ef-\ef^0_r))+\lambda_{(l^0_{r}-[n^\nu]; l^0_{r})} \hat \eo^t_{(l^0_{r}-[n^\nu]; l^0_{r})} |\ef| \\
 \qquad \qquad  +
\sum^{l^0_{r}+[n^\nu]}_{i=l^0_{r}+1}\rho_\tau(\varepsilon_i-b-\XX^t_i(\ef^0_{r+1}))  + \lambda_{(l^0_{r}; l^0_{r}+[n^\nu])} \hat \eo^t_{(l^0_{r}; l^0_{r}+[n^\nu])} |\ef|\}.\end{array}
\]
Applying Lemma \ref{lemma 7aLASSO} for $c_n=c$ on the one of the intervals $(l^0_{r}-[n^\nu]; l^0_{r})$ or $(l^0_{r}; l^0_{r}+[n^\nu])$, (it is the one  where $|b-b^0|+\| \ef-\efo\|_2 > \tilde c >0$), we have $D_4=O_{\eP}(n^\nu)>0$. Then, with probability converging to  1, as $n \rightarrow \infty$, we have
\begin{equation}
\label{ineg}
\inf_{(l_1, \cdots, l_K) \in {\cal L}^c_r(\nu)} \cro{S^*(l_1,\cdots, l_K, l^0_1, \cdots , l^0_{r-1}, l^0_r-[n^\nu], l^0_r+[n^\nu], l^0_{r+1}, \cdots, l^0_K )- S^*_0 }> O_{\eP}(n^\nu).
\end{equation}
Therefore, we have proved that $\eP [ (\hat l^*_1, \cdots, \hat l^*_K) \in {\cal L}^c_r(\nu)] \rightarrow 0$, as $n \rightarrow \infty$.\\
  \textbf{\underline{Step III.}} We now prove that $\hat l^*_r-l^0_r=O_{\eP}(1)$, for each $r=1, \cdots, K$. \\
   Let the set: ${\cal L}(\nu) \equiv \acc{(l_1, \cdots, l_K); |l_t-l^0_t| <[n^\nu], \forall t=1, \cdots, K}$ with $\nu <1/4$. As a consequence of the Step II, for $n$ large, the estimator $(\hat l^*_1, \cdots, \hat l^*_K)$ belongs to ${\cal L}(\nu)$ with a probability tending to 1. We use the reduction to absurdity method, supposing that there exists a change-point estimator at an unbounded distance from the true value. Consider then, for a ${\cal M}_1 >0$ (to be determine later) the following  set: 
\[{\cal L}_r(\nu, {\cal M}_1) \equiv \acc{(l_1, \cdots l_K) \in {\cal L}(\nu); l_r-l^0_r <-{\cal M}_1}.\]
 The case $l_r-l^0_r >{\cal M}_1$ is similar. We shall find a ${\cal M}_1$ such that the probability that the change-point estimator belongs to the set ${\cal L}_r(\nu, {\cal M}_1)$ converges to 0. \\
  Consider two vectors of change-points  $(m_1, \cdots, m_K) \in {\cal L}(\nu)$ and $(l_1, \cdots, l_K) \in {\cal L}_r(\nu,{\cal M}_1)$ such that  $m_t=l_t$ for $t \neq r$ and $m_r=l^0_r$. We have for the difference of the corresponding objective functions:
  \[
 \begin{array}{l} 
  S^*(l_1, \cdots, l_K)- S^*(m_1, \cdots, m_K) =\{ \sum^{l_r}_{i=l_{r-1}+1} [\rho_\tau(Y_i-\hat b^*_{(l_{r-1},l_r)}- \XX^t_i \hat \ef^*_{(l_{r-1},l_r)}) \\
  - \rho_\tau(Y_i-\hat b^*_{(l_{r-1},l^0_r)}- \XX^t_i \hat \ef^*_{(l_{r-1},l^0_r)} )] 
   + \lambda_{(l_{r-1},l_r)} \hat  \eo^t_{(l_{r-1},l_r)} |\hat \ef^*_{(l_{r-1},l_r)} |  -\lambda_{(l_{r-1},l^0_r)} \hat  \eo^t_{(l_{r-1},l^0_r)} |\hat \ef^*_{(l_{r-1},l^0_r)} | \}\\
   +\{\sum^{l^0_r}_{i=l_r+1} [
   \rho_\tau(Y_i-\hat b^*_{(l_{r},l_{r+1})}- \XX^t_i \hat \ef^*_{(l_{r},l_{r+1})})- \rho_\tau(Y_i-\hat b^*_{(l_{r-1},l^0_r)}- \XX^t_i \hat \ef^*_{(l_{r-1},l^0_r)} )] \} \\
  + \{ \sum^{l_{r+1}}_{i=l^0_{r+1}} [ \rho_\tau(Y_i-\hat b^*_{(l_{r},l_{r+1})}- \XX^t_i \hat \ef^*_{(l_{r},l_{r+1})})- \rho_\tau(Y_i-\hat b^*_{(l^0_{r},l_{r+1})}- \XX^t_i \hat \ef^*_{(l^0_{r},l_{r+1})} )] \\
  + \lambda_{(l_r,l_{r+1})} \hat  \eo^t_{(l_r,l_{r+1})} |\hat \ef^*_{(l_r,l_{r+1})} |  -\lambda_{(l^0_r,l_{r+1})} \hat  \eo^t_{(l^0_r,l_{r+1})} |\hat \ef^*_{(l^0_r,l_{r+1})} |
   \} \\
   \equiv \{S_{11}+S_{12} \} +\{ S_{21}\}+\{S_{31}+S_{32}\}.
   \end{array}
\]
    Taking into account  Remark \ref{Remark 1}, is easily obtained that $S_{11}=O_{\eP}(1)$, uniformly in ${\cal M}_1$. Taking into account that $\lambda_{(l_{r-1},l_r)} =O((l_r-l_{r-1})^{1/2})$ together the properties of strong convergence of the  quantile estimator $\hat \ef_{(l_{r-1},l_r)}$ and of adaptive quantile estimator $\hat \ef^*_{(l_{r-1},l_r)}$, we obtain that $S_{12}=o_{\eP}(l^0- l_r)$. Similarly, we have that $S_{31}=O_{\eP}(1)$ and $S_{32}=o_{\eP}(l^0- l_r)$.\\
    It remains to study the most difficult part, that is, $S_{21}$, which can be written\\
$S_{21}=\sum^{l^0_r}_{i=l_r+1}R_i^{(\tau)}(b^0_{r+1},\ef^0_{r+1};b^0_r,\ef^0_r)
 - \sum^{l^0_r}_{i=l_r+1} [   \rho_\tau(\varepsilon_i-\hat b^*_{(l_{r},l_{r+1})}- \XX^t_i (\hat \ef^*_{(l_{r},l_{r+1})}- \ef^0_r)) - \rho_\tau(\varepsilon_i- b^0_{r+1}-\XX^t_i(\ef^0_{r+1}-\ef^0_r))]$ $
   -\sum^{l^0_r}_{i=l_r+1}
  R_i^{(\tau)}(\hat b^*_{(l_{r-1},l^0_r)},\hat \ef^*_{(l_{r-1},l^0_r)} ; b^0_r, \ef^0_r ) $ $
    \equiv S_{211}-S_{212}- S_{213}$.\\
     For $S_{212}$ and $S_{213}$ we use the inequalities 
     \[\left| \frac{\rho_\tau(r_1)- \rho_\tau(r_2)}{r_1-r_2} \right| \leq \max(\tau,1-\tau) < 1\]
      and we obtain that $| S_{212} |$ is smaller than 
\[
  \sum^{l^0_r}_{i=l_r+1} | b^0_{r+1} - \hat b^*_{(l_r,l_{r+1})} +\XX^t_i (\ef^0_{r+1}-\hat \ef^*_{(l_r,l_{r+1})})  | \leq | b^0_{r+1} - \hat b^*_{(l_r,l_{r+1})}|(l^0_r-l_r) +\| \ef^0_{r+1}-\hat \ef^*_{(l_r,l_{r+1})} \|_2 \sum^{l^0_r}_{i=l_r+1}  \|\XX_i\|_2\]
 that is, using Remark \ref{Remark 1}, of order $o_{\eP}(1) $. We obtain analogously $|S_{213}|=o_{\eP}(1)$. For $S_{211}$, combining the relation (\ref{eq_lemma7.3.Scad}) for $c_n=c=\max(|b^0_r-b^0_{r+1}|,\|\ef^0_r - \ef^0_{r+1}\|_2)$ together with Lemma \ref{Lemma 9Bai}  yield that $S_{211} \geq (l^0_r-l_r)\eta \geq {\cal M}_1 \eta$, with probability converging to 1 as ${\cal M}_1 \rightarrow \infty$. \\
 Choosing ${\cal M}_1 >0$ such that $S_{211} \geq \max (|S_{21}|,|S_{11}|,|S_{31}|,|S_{32}|)$ we have that \[\lim_{n \rightarrow \infty} \eP[(\hat l^*_1, \cdots, \hat l^*_K)\in {\cal L}_r(\nu, {\cal M}_1)]=0.\]
  Which proves that  $\hat l^*_r-l^0_r=O_{\eP}(1)$, for each $r=1, \cdots, K$.
 \hspace*{\fill}$\blacksquare$ \\  
 
 \noindent {\bf Proof of Theorem \ref{Theorem oracle}}
 \textit{(i) } This assertion follows immediately from Theorem \ref{Theorem 4.1 ZouYuan}(i)
  and Theorem \ref{Theorem 5aLASSO}.  \\
\hh  \textit{(ii)} 
  \begin{figure}[h!]
\begin{minipage}[b]{0.50\linewidth}
\includegraphics[scale=0.38]{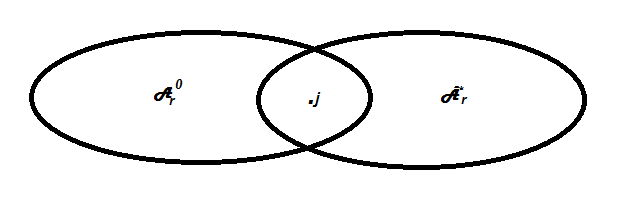} 
\caption{}
\label{Figure 2}
\end{minipage}
\begin{minipage}[b]{0.50\linewidth}
\includegraphics[scale=0.38]{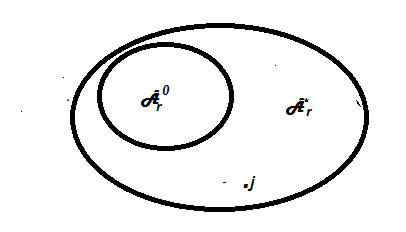}
\caption{}
\label{Figure 3}
\end{minipage}
\end{figure} 
 The consistency of the variable selection in a model with one phase, property established by Theorem \ref{Theorem 4.1 ZouYuan}(ii), implies that 
 \begin{equation}
 \label{spar1}
 \lim_{n \rightarrow \infty} \eP \cro{\hat {\cal A}^0_{n,r} = {\cal A}^0_r}=1.
 \end{equation}
\hh   It remains to prove that $\lim_{n \rightarrow \infty} \eP \cro{\hat {\cal A}^0_{n,r} =\hat {\cal A}^*_{n,r}}=1$.
 The general case is considered  in the early, that is presented in Figure \ref{Figure 2}. \\
\hh   If $j \in \hat {\cal A}^0_{n,r}$, thus, using (\ref{spar1}), we have that $j \in {\cal A}^0_r$ with probability tending to  1, which implies  that $\phi^0_{r,j} \neq 0$. Moreover, using the result proved in the previous question (i), and the fact that  the  adaptive LASSO quantile estimator  $\hat \phi^*_{(\hat l^*_{r-1};\hat l^*_r),j}$ of the $j$th component of the regression parameter $\ef^0_r$, calculated between the corresponding adaptive LASSO quantile estimators of the change-points,  is asymptotically normal, we obtain that $\hat \phi^*_{(\hat l^*_{r-1};\hat l^*_r);j} \overset{\eP} {\underset{n \rightarrow \infty}{\longrightarrow}} \phi^0_{r,j} \neq 0$. Then,  
\[\lim_{n \rightarrow \infty} \eP \cro{\hat \phi^*_{(\hat l^*_{r-1};\hat l^*_r);j} \neq 0 }=1,  \textrm{ i.e. } \eP[ j \in \hat {\cal A}^*_{n,r}] \rightarrow 1.\]
 Thus $\eP[ \hat {\cal A}^0_{n,r} \subseteq \hat {\cal A}^*_{n,r} ] \rightarrow 1$.\\
\hh There, remains now  the most difficult part: to prove that, if the index $j \not \in \hat {\cal A}^0_r $, then $j \not \in \hat {\cal A}^*_{n,r}$ (in fact, in view of (\ref{spar1}), we must show that if the true component is zero, the corresponding estimators don't converge to 0 and  for a fixed $n$, the component $\hat \phi^*_{(\hat l^*_{r-1};\hat l^*_r);j} \neq 0$). We will then calculate $\eP [j \not \in {\cal A}^0_r, j \in \hat {\cal A}^*_{n,r} ]$ (see Figure \ref{Figure 3}). Since  $j \in \hat {\cal A}^*_{n,r}$, by KKT optimality condition, Proposition \ref{KKT}(i), we have
 \[
  \tau \sum^{\hat l^*_r}_{i=\hat l^*_{r-1}+1} X_{ij}- \sum^{\hat l^*_r}_{i=\hat l^*_{r-1}+1} X_{ij} \e1_{Y_i < \XX^t_i \hat \ef^*_{(\hat l^*_{r-1};\hat l^*_r)}} = \lambda_{(\hat l^*_{r-1};\hat l^*_r)} \hat \omega_{(\hat l^*_{r-1};\hat l^*_r),j} sgn(\hat \phi^*_{(\hat l^*_{r-1};\hat l^*_r),j}).
 \]
 Then 
 \begin{equation}
 \label{eq17}
 \lambda_{(\hat l^*_{r-1};\hat l^*_r)} \hat \omega_{(\hat l^*_{r-1};\hat l^*_r),j} = \left| \tau \sum^{\hat l^*_r}_{i=\hat l^*_{r-1}+1} X_{ij}- \sum^{\hat l^*_r}_{i=\hat l^*_{r-1}+1} X_{ij} \e1_{Y_i < \XX^t_i \hat \ef^*_{(\hat l^*_{r-1};\hat l^*_r)}} \right| 
 < 2 \sum^{\hat l^*_r}_{i=\hat l^*_{r-1}+1} | X_{ij} |. 
 \end{equation}
 On the other hand
 \begin{equation}
 \label{eq17bis}
 \frac{\lambda_{(\hat l^*_{r-1};\hat l^*_r)} \hat \omega_{(\hat l^*_{r-1};\hat l^*_r),j} }{\hat l^*_r -\hat l^*_{r-1} } =  \frac{\lambda_{(\hat l^*_{r-1};\hat l^*_r)}}{|(\hat l^*_r -\hat l^*_{r-1})^{1/2} \hat \phi_{(\hat l^*_{r-1};\hat l^*_r),j} |^g} \cdot \frac{(\hat l^*_r -\hat l^*_{r-1})^{g/2}}{\hat l^*_r -\hat l^*_{r-1}},
 \end{equation}  
 where  $\hat \phi_{(\hat l^*_{r-1};\hat l^*_r),j}$ is the quantile estimator of $\phi^0_{r,j}$. Since  $j  \in {\cal A}^0_r$, we have that  $\phi^0_{r,j} \neq 0$. On the other hand, taking into account Theorem \ref{Theorem 5aLASSO} and the fact that the quantile estimator $\hat \phi_{(\hat l^*_{r-1};\hat l^*_r),j}$ is strongly consistent and asymptotically normal, we have that for all $\epsilon >0$, there exists $\eta_\epsilon >0$ such that 
 \begin{equation}
 \label{eq16}
 \eP[\pth{(\hat l^*_r -\hat l^*_{r-1})^{1/2} | \hat \phi_{(\hat l^*_{r-1};\hat l^*_r),j} | }^{-1} >\eta_\epsilon ] >1 -\epsilon . 
\end{equation}  
Since $(l_r-l_{r-1})^{g/2-1} \lambda_{(l_{r-1}, l_r)} \rightarrow \infty$, together with the relation (\ref{eq16}), we obtain that (\ref{eq17bis}) converges to infinity with probability converging to 1 as $n \rightarrow \infty$. On the other hand, an application of  Cauchy-Schwarz's inequality yields that 
\[
(\hat l^*_r -\hat l^*_{r-1})^{-1} \sum^{\hat l^*_r}_{i=\hat l^*_{r-1}+1} | X_{ij} | \leq \pth{(\hat l^*_r -\hat l^*_{r-1})^{-1} \sum^{\hat l^*_r}_{i=\hat l^*_{r-1}+1}  X^2_{ij} }^{1/2},\] 
wich is bounded with probability converging to 1, from assumption (A1). Taking into account (\ref{eq17}), we obtain that (\ref{eq17bis}) is bounded. Contradiction. Then 
\[\eP [j \not \in {\cal A}^0_r, j \in \hat {\cal A}^*_{n,r} ] \rightarrow 0, \textrm{ as } n \rightarrow \infty.
\]
 Which implies that $\eP[\hat {\cal A}^*_{n,r} \subseteq   {\cal A}^0_r ] \rightarrow 1 $ and in view of relation (\ref{spar1}), we have that $\eP[\hat {\cal A}^*_{n,r} \subseteq \hat  {\cal A}^0_{n,r} ] \rightarrow 1 $, as $n \rightarrow \infty$.
   \hspace*{\fill}$\blacksquare$ \\
  
  \noindent {\bf Proof of Theorem  \ref{Theorem 4aLASSO}}
For $K^0$ the true number of changes, let us define the objective function calculated for the true values of the parameters. Only the weights are estimated:
\begin{equation}
\label{def_S0*}
{\cal S}_0 \equiv \sum^n_{i=1} \sum^{K^0+1}_{r=1} \rho_\tau(\varepsilon_i - b^0_r \e1_{l^0_{r-1} \leq i < l^0_r}) +\sum^{K^0+1}_{r=1} \lambda_{(l^0_{r-1};l^0_r)} \hat \eo^t_{(l^0_{r-1};l^0_r)} |\ef^0_r|, 
\end{equation}   
with $\hat \eo^t_{(l^0_{r-1};l^0_r)} = |\hat \ef_{(l^0_{r-1};l^0_r)} |^{-g}$ , calculated  on the basis of the  quantile estimator, on the observations between $l^0_{r-1}$ et $l^0_r$.  \\
\hh We will first study the behavior of  $\hat s^*_{K^0}$ for the true number of phases. 
Using Lemma \ref{lemma 6aLASSO}, for $\alpha > 1/2$, we have for the difference between the objective function  $S^*$ calculated for the adaptive LASSO change-points estimators and the sum calculated for the true values:
\begin{equation}
\label{eq22}
 S^*(\hat l^*_{1,K^0}, \cdots , \hat l^*_{K^0,K^0})- {\cal S}_0 = \min (o_{\eP}(n^{(1+g)/2}), O_{\eP}(n^\alpha) ).
\end{equation}
Then 
\[\hat s^*_{K^0} = n^{-1} S^*(\hat l^*_{1,K^0}, \cdots , \hat l^*_{K^0,K^0})- n^{-1} {\cal S}_0 +n^{-1} {\cal S}_0=\min (o_{\eP}(n^{(g-1)/2}), O_{\eP}(n^{\alpha-1/2}) )+n^{-1} {\cal S}_0.\]
 By the weak law of large numbers, using assumption (A4) and the independence of $(\varepsilon_i)$, we have \[n^{-1} {\cal S}_0 \overset{\eP} {\underset{n \rightarrow \infty}{\longrightarrow}} \sum^{K^0+1}_{r=1} \eE[\rho_\tau(\varepsilon-b^0_r) ].\]
  Thus, 
\begin{equation}
\label{eq23}
\hat s^*_{K^0} \overset{\eP} {\underset{n \rightarrow \infty}{\longrightarrow}} \sum^{K^0+1}_{r=1} \eE[\rho_\tau(\varepsilon-b^0_r) ] > 0.
\end{equation}
\hh Now,  we show that
\begin{equation}
\label{eq20}
\eP[ \hat K^*_n < K^0] \rightarrow 0, \qquad \textrm{ as } n \rightarrow \infty.
\end{equation}
In order to prove (\ref{eq20}), consider $K$ any change-point number , with  $K < K^0$. Then
\[
B(K)-B(K^0)= n \log \pth{ 1+ \frac{\hat s^*_K -\hat s^*_{K^0}}{\hat s^*_{K^0}}}+B_n[ G(K,p) - G(K^0,p) ]
.\]
 Two cases are possible concerning $(\hat s^*_K -\hat s^*_{K^0})/\hat s^*_{K^0}$.  
\begin{itemize}
\item If $ (\hat s^*_K -\hat s^*_{K^0})/\hat s^*_{K^0}$ is great ($ \geq C  > 0$), then, since $\hat s^*_{K^0} > \epsilon_2 > 0$, we have that there exists $\epsilon_1 > 0$ such that $\hat s^*_K -\hat s^*_{K^0} > \epsilon_1 > 0$ for any $n$ large enough. Then 
\[\eP[\argmin_K  \hat s^*_K <K^0] \rightarrow 0, \textrm{ for } n \rightarrow \infty \]
 and the relation (\ref{eq20}) follows.
\item If $(\hat s^*_K -\hat s^*_{K^0})/\hat s^*_{K^0}=o_{\eP}(1)$, then, using the fact that for  $x$ close to 0 we have $\log(1+x) \simeq x$, we have
\begin{equation}
\label{eq21}
B(K)-B(K^0)=n \frac{\hat s^*_K -\hat s^*_{K^0}}{\hat s^*_{K^0}}(1+o_{\eP}(1))+B_n[ G(K,p) - G(K^0,p) ].
\end{equation}
We will study in this case the first term of the right side of the relation (\ref{eq21}). Recall that the constant $a \in (1/2,1)$ is that of the assumption (A3). Similarly of  the relation (\ref{ineg}), since between  any two  consecutive  change-points $l_{r-1}$ and $l_r$ there are at least $[n^a]$ observations and since when $K < K^0$ there is at least a non estimated true change-point, we obtain that
\begin{equation}
\label{eq24}
S^*(\hat l^*_{1,K}, \cdots , \hat l^*_{K,K}) -{\cal S}_0 > C n^a.
\end{equation} 
Then, $n(\hat s^*_K -\hat s^*_{K^0}) /\hat s^*_{K^0}= [S^*(\hat l^*_{1,K}, \cdots , \hat l^*_{K,K}) -{\cal S}_0 - S^*(\hat l^*_{1,K^0}, \cdots , \hat l^*_{K^0,K^0}) +{\cal S}_0  ]/\hat s^*_{K^0}$. Using  relations (\ref{eq22}), (\ref{eq23}), for $a >\alpha>1/2$, we obtain 
\begin{equation}
\label{eq25}
n(\hat s^*_K -\hat s^*_{K^0}) /\hat s^*_{K^0}> C [ O_{\eP}(n^a)-O_{\eP}(n^\alpha)]= O_{\eP}(n^a). 
\end{equation}
Using the relations (\ref{eq21}), (\ref{eq25}),  the fact that $B_n=o(n^a)$ and since the function $G$ is increasing in $K$, we obtain that for $K < K^0$ we have  $B(K)-B(K^0)>  O_{\eP}(n^a)- o_{\eP}(n^a) \rightarrow \infty$. Thus, the relation (\ref{eq20}) follows.
\end{itemize}
\hh We finally consider the case  $K > K^0$, cases wherein, given the definition of $S^*$ and of the change-points estimators, we have
\begin{equation}
\label{eq26}
\begin{array}{ll}
{\cal S}_0 \geq S^*(l^0_1, \cdots, l^0_{K^0}) & \geq S^*(\hat l^*_{1,K^0}, \cdots , \hat l^*_{K^0,K^0}) \geq  S^*(\hat l^*_{1,K}, \cdots , \hat l^*_{K,K}) \\
 & \geq S^*(\hat l^*_{1,K}, \cdots , \hat l^*_{K,K},l^0_1, \cdots, l^0_{K^0} )
 \end{array}
\end{equation}
Then, by similar calculations as for the inequality (\ref{ineg}) of Theorem \ref{Theorem 5aLASSO}, we have that for the last term of (\ref{eq26}): 
\[S^*(\hat l^*_{1,K}, \cdots , \hat l^*_{K,K},l^0_1, \cdots, l^0_{K^0} )-{\cal S}_0 > O_{\eP}(n ^\nu), \qquad  \textrm{ for } \nu <1/4.\]
 Thus $0 \geq \hat s^*_{K^0} - \hat s^*_{K}= O_{\eP}(n^{\nu-1})$. Since $G(K,p)$ increases in $K$, together  $n^\alpha \ll B_n \ll n^a$ and the relation (\ref{eq23}), we have that 
 \[n \log \hat s^*_{K^0} - n \log \hat s^*_{K} = n (\hat s^*_{K^0} - \hat s^*_{K})(1+o_{\eP}(1))=O_{\eP}(n^\nu) .\]
  In fact, it's the penalty that takes over in this case  $K > K^0$. In this circumstance, we have  $G(K, p) \geq G(K^0,p)$ and then  $B_n[ G(K, p) -  G(K^0,p)]\geq C B_n > O(n^{1/2})$. Hence, we have for the difference between the values of the two  criteria $B(K)-B(K^0) \gg - O_{\eP}(n^\nu)+O(n^{1/2})= O(n^{1/2}) \rightarrow \infty $ as $n\rightarrow \infty$. Then, for large enough $n$, $n \log \hat s^*_{K} +G(K, p) > n \log \hat s^*_{K^0}   +G(K^0,p)$ implying 
  \begin{equation}
  \label{AFC1}
  \eP[\hat K_n >K^0] \rightarrow 0, \qquad \textrm{  as } n\rightarrow \infty.
  \end{equation}
  The Theorem follows from relations (\ref{eq20}) and (\ref{AFC1}). 
  \hspace*{\fill}$\blacksquare$ \\ 
  
\subsection{Lemma  proofs}
 \noindent {\bf Proof of Lemma \ref{Lemma 8Bai}}
 Consider the notations $\tilde \ef \equiv (b-b^0, \ef-\efo)$,  $w_i=(1,\XX^t_i)$ and afterward we apply Lemma 8 of \cite{Bai:98}.
  \hspace*{\fill}$\blacksquare$ \\ 
  
\noindent {\bf Proof of Lemma \ref{lo}}\\
\hh  \underline{Case I:} there exists a deterministic  sequence $(d_n)$ that tends to infinity,  as $n \rightarrow \infty$, such that $(k-l) \geq d_n $.\\
 By the  definition of $\hat \eo_{(l;k)} = | \hat \ef_{(l;k)}|^{-g}$, if there exists a  component $j$ of the quantile estimator such that $\hat \phi_{(l;k),j} \rightarrow 0$ (or a subsequence),  then, since $\eP[\hat \phi_{(l;k),j} =0] =0$ and $\hat \omega_{(l;k),j}=|\hat \phi_{(l,k),j}|^{-g} $, using the asymptotic normality of  $\hat \ef_{(l;k)}$ we have that $\hat \phi_{(l,k),j}=O_{\eP}(k-l)^{-1/2}=O_{\eP}(n^{-1/2})$. Thus, $\lambda_{(l;k)} \hat \omega_{(l;k),j}=o(n^{1/2}) O_{\eP}(n^{g/2})=o_{\eP}(n^{(1+g)/2})$. If the $j$-th component is such that  $\hat \phi_{(l;k),j} \geq c >0 $ then  $\lambda_{(l;k)} \hat \omega_{(l;k),j}=o_{\eP}(n^{1/2})$.\\
\hh \underline{Case II. } If $(k-l)$ does not converge  with $n$ to infinity, then, we can extract a bounded subsequence. 
Let us suppose that $0 \leq l < k \leq c$. Since $\eP[\hat \phi_{(l;k),j} =0] =0$, we have in this case $\lambda_{(l;k)} \hat \omega_{(l;k),j}=o_{\eP}(n^{1/2})$.
  \hspace*{\fill}$\blacksquare$ \\

 \noindent {\bf Proof of  Lemma \ref{lemma 6aLASSO}}
 Let us denote the difference $k-l$ by $d_n$. 
Taking into account Lemma \ref{lo}, then  Lemma \ref{lemma 6aLASSO} is proven if an equivalent result for  $R_i^{(\tau)}$ is showed:
 \begin{equation}
 \label{eq8}
 \sup_{0 \leq l < k \leq n} \left| \inf_{b, \ef} \sum^k_{i=l+1} R_i^{(\tau)}(b,\ef;b^0,\efo) \right| =O_{\eP}(n^\alpha), \qquad \alpha >1/2.
 \end{equation}
Let us consider the random  processes   
\[H^{(\tau)}_i(b,\ef; b^0,\efo) \equiv R^{(\tau)}_i(b,\ef; b^0,\efo)+ [\tau \e1_{\varepsilon_i > b^0 }- (1-\tau) \e1_{\varepsilon_i \leq b^0} ] \XX^t_i(\ef-\efo).\]
 Obviously $\eE[\tau \e1_{\varepsilon_i > b^0 }- (1-\tau) \e1_{\varepsilon_i \leq b^0} ]=(1-\tau) \eP[\varepsilon_i \leq b^0] - \tau \eP[\varepsilon_i > b^0]=0$.  \\
 On the other hand, we have, for all $b_1, b_2 \in {\cal B}$ and $\ef_1, \ef_2 \in \Gamma$, that \\ $H^{(\tau)}_i(b_1,\ef_1; b^0,\efo)-H^{(\tau)}_i(b_2,\ef_2; b^0,\efo)=\rho_\tau(\varepsilon_i-b_1-\XX^t_i(\ef_1-\efo))- \rho_\tau(\varepsilon_i-b_2-\XX^t_i(\ef_2-\efo))+ [\tau \e1_{\varepsilon_i > b^0 }- (1-\tau) \e1_{\varepsilon_i \leq b^0} ] \XX^t_i(\ef_1-\ef_2)$. But, generally, for all $r_1, r_2 \in \R$ we have that $|\rho_\tau(r_1) - \rho_\tau(r_2)| < |r_1-r_2|$. Thus
 \[|\rho_\tau(\varepsilon_i-b_1-\XX^t_i(\ef_1-\efo))- \rho_\tau(\varepsilon_i-b_2-\XX^t_i(\ef_2-\efo))| < |\XX^t_i(\ef_1-\ef_2)+b_1-b_2 |.\]
  If the parameters $b_1, b_2, \ef_1, \ef_2  $ are such that  $\|\ef_1- \ef_2\|_2 \leq C n^{-1/2}$, $|b_1-b_2| \leq C n^{-1/2}$, since $\eE[\tau \e1_{\varepsilon_i > b^0 }- (1-\tau) \e1_{\varepsilon_i \leq b^0} ]=0$, we have, using  assumption  (A1), that
  \begin{equation}
 \label{eq5}
 \begin{array}{l}
 \sum^n_{i=1} \left\{\rho_\tau(\varepsilon_i-b_1-\XX^t_i(\ef_1-\efo))- \rho_\tau(\varepsilon_i-b_2-\XX^t_i(\ef_2-\efo)) \right. \\
\quad   \left.-\eE[\rho_\tau(\varepsilon_i-b_1-\XX^t_i(\ef_1-\efo))] 
 + \eE[\rho_\tau(\varepsilon_i-b_2-\XX^t_i(\ef_2-\efo))] \right\} \\
 \quad \leq C \sum^n_{i=1} |b_1- b_2| +C \sum^n_{i=1} \| \XX_i\|_2 \cdot \| \ef_2-\ef_1 \|_2 =O_{\eP}(n^{1/2}).
 \end{array}
 \end{equation}
On the other hand, by Proposition \ref{proposition 2.1.SCAD}, we have that $\eE[R^{(\tau)}_i(b,\ef; b^0,\efo)] \geq 0$, which means, taking into account the fact that $R^{(\tau)}_i(b^0,\efo; b^0,\efo)=0$, that, for $k-l=d_n \rightarrow \infty$, we have 
\[
 0 \geq \inf_{b, \ef}  \sum^{k}_{i=l+1} R_{i}^{(\tau)}(b,\ef;b^0,\efo) \geq \inf_{b, \ef}  \sum^{k}_{i=l+1}  [ R_{i}^{(\tau)}(b,\ef;b^0,\efo) -\eE[R_{i}^{(\tau)}(b,\ef;b^0,\efo)]].\]
Thus 
 \[
 | \inf_{b, \ef} \sum^{k}_{i=l+1} R_{i}^{(\tau)}(b,\ef;b^0,\efo)| \leq \sup_{b, \ef} | \sum^{k}_{i=l+1} [ R_{i}^{(\tau)}(b,\ef;b^0,\efo) -\eE[R_{i}^{(\tau)}(b,\ef;b^0,\efo)]] |.\]
Then, for  $d_n \rightarrow \infty $, we have
 \begin{equation}
 \label{AA}
 \sup_{0 \leq l < k \leq n, k-l=d_n} \left| \inf_{b, \ef} \sum^k_{i=l+1} R_i^{(\tau)}(b,\ef;b^0,\efo) \right| \leq 2 \sup_{d_n \leq k \leq n} \zeta_k ,
 \end{equation}
with $\zeta_k \equiv \sup_{b,\ef} | \sum_{i=1}^k [R_i^{(\tau)}(b,\ef;b^0,\efo)- \eE[R_i^{(\tau)}(b,\ef;b^0,\efo) ]] | $. Thus, by Proposition \ref{proposition 2.1.SCAD}, $\{ \zeta_k,{\cal F}_k\}$ is a submartingale, where ${\cal F}_k \equiv \sigma-field \{ \varepsilon_1, \cdots, \varepsilon_k\}$. Which  means, by Doob's inequality, for $\alpha > 1/2$, that:
 \begin{equation}
 \label{eq6}
 \eP \cro{\sup_{d_n \leq k \leq n} \zeta_k > n^\alpha } \leq \eP \cro{\sup_{1 \leq k \leq n}  \zeta_k > n^\alpha } \leq n^{- \alpha m} C_m \eE[\zeta_n^m ],
 \end{equation}
for some $C_m >0$, and $m>1$ determined later. \\
 Now divide the parameter set ${\cal B} \times \Gamma$ in $n^{(1+p)/2}$ cells, such that the diameter of each cell is less than $n^{-1/2}$. Thus, for  $(b_1,\ef_1), (b_2,\ef_2)$ in the same cell, we have:
\begin{equation}
 \label{eq7}
 \begin{array}{l}
 \sum^n_{i=1} \cro{  R_{i}^{(\tau)}(b_1,\ef_1;b^0,\efo) -\eE[R_{i}^{(\tau)}(b_1,\ef_1;b^0,\efo)]- R_{i}^{(\tau)}(b_2,\ef_2;b^0,\efo) + \eE[R_{i}^{(\tau)}(b_2,\ef_2;b^0,\efo)]} \\
 \qquad  \leq \left|\sum^n_{i=1} \left[ H_{i}^{(\tau)}(b_1,\ef_1;b^0,\efo) -\eE[H_{i}^{(\tau)}(b_1,\ef_1;b^0,\efo)]- H_{i}^{(\tau)}(b_2,\ef_2;b^0,\efo) \right. \right.  \\
  \qquad \qquad \qquad +\left. \left.  \eE[H_{i}^{(\tau)}(b_2,\ef_2;b^0,\efo)] \right] \right| + \left| \sum^n_{i=1}  D_i \XX^t_i (\ef_1-\ef_2) \right| \leq C n^{1/2},
 \end{array}
  \end{equation}
where the last inequality follows from (\ref{eq5}). Let be now $(b_j,\ef_j)$ in the  $j$th cell, $j=1, \cdots, c_p n^{(1+p)/2}$. Then,  as in the paper of \cite{Bai:98}, Lemma 3, we have 
\[
 \zeta^m_n \leq C \sup_j | \sum^n_{i=1} [R^{(\tau)}_i(b_j,\ef_j ; b^0, \efo) - \eE[R^{(\tau)}_i(b_j,\ef_j ; b^0, \efo)] ]  |^m+C n^{m/2}.\]  
  Since $R^{(\tau)}_i(b_j,\ef_j ; b^0, \efo) - \eE[R^{(\tau)}_i(b_j,\ef_j ; b^0, \efo)]$ is a bounded martingale difference, for each fixed $j \in\{ 1, \cdots, c_p n^{(1+p)/2} \} $, we have  
  \[\eE [| \sum^n_{i=1} (R^{(\tau)}_i(b_j,\ef_j ; b^0, \efo) - \eE[R^{(\tau)}_i(b_j,\ef_j ; b^0, \efo)] ) |^m] \leq C n^{m/2}.\]
   Then $\eE[\zeta^m_n] \leq O_{\eP}(n^{\frac{1+p}{2}}) n^{m/2}+O_{\eP}(n^{m/2})=O_{\eP}(n^{\frac{m+p+1}{2}})$.\\
   Thus, choosing $m$ such that $m >\frac{p+1}{2 \alpha +1}$, the right-hand side of (\ref{eq6}) converges to 0. Furthermore, taking into account  (\ref{AA}), (\ref{eq6}) and (\ref{eq7}) we obtain  the claim (\ref{eq8}).  
 \hspace*{\fill}$\blacksquare$ \\
 
    \noindent {\bf Proof of  Lemma  \ref{Lemma 4Bai}} 
Let us consider the set  $\Omega_n$ written as an union of subsets  $\Omega_n=\bigcup^{C_p n^{(p+1)/2}}_{j=1} \mathcal{C}_j^n$, with $C_p $ a bounded positive constant,  and 
\[{\cal C}^n_j \equiv \{(b,\ef) \in \Omega_n; |b-b'|+\|\ef- \ef'\|_2 \leq c_n n^{-1/2},   \textrm{ for all } (b',\ef') \in {\cal C}_j^n \}.\]
 For $(b_1,\ef_1),(b_2,\ef_2) \in {\cal C}_j^n$, as in the proof of Lemma \ref{lemma 6aLASSO}, relation (\ref{eq7}), we have that
$(n c^2_n)^{-1}|{\cal R}_n^{(\tau)}(b_1,\ef_1;b^0,\efo)- \eE[{\cal R}_n^{(\tau)}(b_2,\ef_2;b^0,\efo)] -{\cal R}_n^{(\tau)}(b_2,\ef_2;b^0,\efo)+ \eE[{\cal R}_n^{(\tau)}(b_2,\ef_2;b^0,\efo)]|
$
$ \leq  C (|b_1-b_2|-\|\ef_1- \ef_2\|_2)c^{-2}_n \leq C n^{-1/2}c^{-1}_n \rightarrow 0$, as $n \rightarrow \infty$. Then, for all $(b_1,\ef_1),(b_2,\ef_2) \in {\cal C}_j^n$,
\begin{equation}
\label{AFC2}
\begin{array}{l}
\lim_{n \rightarrow \infty} \frac{1}{n c^2_n} \left|{\cal R}_n^{(\tau)}(b_1,\ef_1;b^0,\efo)- \eE[{\cal R}_n^{(\tau)}(b_2,\ef_2;b^0,\efo)] -{\cal R}_n^{(\tau)}(b_2,\ef_2;b^0,\efo) \right. \\
\qquad \qquad \qquad \qquad + \left. \eE[{\cal R}_n^{(\tau)}(b_2,\ef_2;b^0,\efo)] \right| =0.
\end{array}
\end{equation}
\hh For $(b_j,\ef_j)\in {\cal C}_j^n$, for all  $j=1, \cdots, C_p n^{(1+p)/2}$, we have that the probability  
\begin{equation}
\label{eq10}
\begin{array}{l}
\eP[\sup_{j}|(n c^2_n)^{-1} [ {\cal R}_n^{(\tau)}(b_j,\ef_j;b^0,\efo)- \eE[{\cal R}_n^{(\tau)}(b_j,\ef_j;b^0,\efo)] ]  | >\epsilon ]  \\
\qquad \leq \sum^{C_p n^{(1+p)/2}}_{j=1} \eP \cro{ \left|{\cal R}_n^{(\tau)}(b_j,\ef_j;b^0,\efo)- \eE[{\cal R}_n^{(\tau)}(b_j,\ef_j;b^0,\efo)] \right| > nc^2_n \epsilon }.
\end{array}
\end{equation}
But, by assumption (A1), we have that 
\[R_i^{(\tau)}(b_j,\ef_j;b^0,\efo) -\eE[R_i^{(\tau)}(b_j,\ef_j;b^0,\efo) ] \leq C [ |b_j-b^0|+ \| \ef_j-\efo\|_2] \leq C c_n.\]
 Then  $Var [ R_i^{(\tau)}(b_j,\ef_j;b^0,\efo)] \leq C^2 c^2_n$ uniformly in $b_j$ and $\ef_j$. We apply  Lemma  \ref{l 1},   for $\beta=C c_n$, $V=C^2 n c_n^2$, $s=1/2$, $z=n c_n^2 \epsilon$ and we obtain 
 \begin{equation}
 \label{AFC3}
 \eP [| {\cal R}_n^{(\tau)}(b_j,\ef_j;b^0,\efo)- \eE[{\cal R}_n^{(\tau)}(b_j,\ef_j;b^0,\efo)] | > n c_n^2 \epsilon ] \leq 2 \exp(-\epsilon^2 nc^2_n C ).
 \end{equation}
The statement of  Lemma \ref{l 1} is given at the beginning of Section 6. 
Relation (\ref{AFC3}) and the fact that $n c^2_n/ \log n \rightarrow \infty$ imply that   the right-hand side of (\ref{eq10}) is bounded by  $2 C_p n^{(1+p)/2}\exp(-\epsilon^2 nc^2_n C )$ which is smaller than $\exp(- \epsilon^2 n c^2_n C/2 )$, for any large enough $n$.\\
Then, for all $(b_j,\ef_j)\in {\cal C}_j^n$, for all  $j=1, \cdots, C_p n^{(1+p)/2}$, we have 
\begin{equation}
\label{AFC4}
 \lim_{n \rightarrow \infty} \eP \left[\sup_{j} \left|(n c^2_n)^{-1} [ {\cal R}_n^{(\tau)}(b_j,\ef_j;b^0,\efo)- \eE[{\cal R}_n^{(\tau)}(b_j,\ef_j;b^0,\efo)] \right]  \right| >\epsilon ] =0 .
\end{equation}
The Lemma follows from (\ref{AFC2}) and (\ref{AFC4}). 
  \hspace*{\fill}$\blacksquare$ \\
 
  \noindent {\bf Proof of  Lemma \ref{lemma 7aLASSO}}
By the proof of   Proposition \ref{proposition 2.1.SCAD}, using  assumption (A2), we have
  \[
 \eE[{\cal R}_n^{(\tau)}(b^0+u_0/\sqrt{n},\efo+\uu/\sqrt{n};b^0,\efo)]=\frac{1}{2 n} f(b^0) (u_0,\uu^t) \left[
   \begin{array}{cc}
   n & \textbf{0}\\
   0 & \sum^n_{i=1}\XX_i \XX^t_i
   \end{array}
   \right] (u_0,\uu^t)^t (1+o(1)),
  \]
 with  $u_0 \in \R$ and $\uu \in \R^p$ in open sets. 
 For a sequence $(c_n)$ converging to zero but with a slower rate than $n^{-1/2}$, we have
  \[
  \eE \cro{\int^{(u_0+\XX^t_i \uu)c_n}_0 [\e1_{\varepsilon_i \leq b^0+t} -\e1_{\varepsilon_i \leq b^0} ] dt} = \int^{u_0+\XX^t_i \uu}_0 c_n [F(b^0+c_n v)- F(b^0)]dv.\] 
Thus
  \[
  \eE [ {\cal R}^{(\tau)}_n((b^0+u_0)c_n, (\ef^0+\uu)c_n;b^0,\efo)] = c_n \sum^n_{i=1} \int^{u_0+\XX^t_i \uu}_0  [F(b^0+c_n v)- F(b^0)]dv 
  \]
  \[
 \qquad \qquad \qquad \qquad  =c^2_n \frac{f(b^0)}{2 }  (u_0,\uu^t) \left[
   \begin{array}{cc}
   n & \textbf{0}\\
   0 & \sum^n_{i=1}\XX_i \XX^t_i
   \end{array}
   \right] (u_0,\uu^t)^t (1+o(1)).
  \]
  By Lemma \ref{Lemma 4Bai}, using the Borel-Cantelli lemma, we have that for any $\epsilon_2>0$
\begin{equation}
\label{eq9}
\limsup_{n \rightarrow \infty} \pth{\sup_{(b,\ef) \in \Omega_n}\left| \frac{1}{n c^2_n } \cro{{\cal R}_n^{(\tau)}(b,\ef;b^0,\efo) - \eE[{\cal R}_n^{(\tau)}(b,\ef;b^0,\efo) ]} \right|} \leq \epsilon_2, \quad a.s.
\end{equation}
But, the function $R_i^{(\tau)}(b,\ef;b^0,\efo)$ is  convex in $b$ and $\ef$, therefore also its sum ${\cal R}_n^{(\tau)}$. Since $ {\cal R}_n^{(\tau)}(b^0,\ef^0;b^0,\efo) =0$, using assumption (A1) and relation (\ref{eq9}), by a similar argument as in Remark 4 of \cite{Bai:98}, we obtain that the infimum of $(b,\ef) \in \Omega^c_n$ is attained on the boundary: $|b-b^0|=c_n$ and $\| \ef-\efo\|_2=c_n$.\\
\hh Then, let us consider the scalar  $u_0$ and the $p$-vector $\uu$, each of euclidean norm 1: $|u_0|=1$, $\| \uu\|_2=1$.  Using  assumption  (A1), we have that:
\begin{equation} 
\label{er}
 \eE \cro{ {\cal R}^{(\tau)}_n((b^0+u_0)c_n, (\ef^0+\uu)c_n;b^0,\efo) ] dt} =n c^2_n \frac{f(b^0)}{2 } (C+o(1)),
\end{equation}
with $C>0$. \\
Then, considering  $\epsilon_2=C f(b^0)/4 $ in (\ref{eq9}), taking into account the relation (\ref{er}), we obtain that
 there exists $\epsilon_3 >0$ (we can consider by example $\epsilon_3=\epsilon_2$) such that, with the probability 1:
\begin{equation}
\label{eq_lemma7.3.Scad}
\liminf_{n \rightarrow \infty} \pth{ \inf_{(b,\ef) \in \Omega^c_n} \frac{1}{n c^2_n}  {\cal R}^{(\tau)}_n (b,\ef;b^0,\efo) } \geq \epsilon_3,
\end{equation}
knowing that $\inf_{(b,\ef) \in \Omega^c_n} $  is in fact the infimum on the boundary of the set $\Omega^c_n $.\\
To finish the lemma, we take into account that:
\[
 \inf_{(b,\ef) \in \Omega^c_n} \frac{1}{n c^2_n} \sum^n_{i=1} R^{(\tau,\lambda)}_i (b,\ef;b^0,\efo) \qquad \qquad \qquad \qquad \qquad \qquad \qquad \qquad \]
 \[ \geq \inf_{(b,\ef) \in \Omega^c_n} \frac{1}{n c^2_n}  {\cal R}^{(\tau)}_n (b,\ef;b^0,\efo) - \frac{\lambda_n}{n c^2_n}  \sup_{(b,\ef) \in \Omega^c_n}  (| \hat \eo^t_{(0;n)}| \cdot \left| |\ef|-|\efo| \right|)
\]
\begin{equation}
\label{ep}
\qquad \qquad \qquad \qquad \geq  \inf_{(b,\ef) \in \Omega^c_n} \frac{1}{n c^2_n}  {\cal R}^{(\tau)}_n (b,\ef;b^0,\efo) - \frac{\lambda_n}{n c^2_n} \| \hat \eo^t_{(0;n)}\|_1 C_+, 
\end{equation}
with $C_+$ a positive constant. The first term of the right-hand side of  (\ref{ep}) is greater than $\epsilon_3$ by the relation (\ref{eq_lemma7.3.Scad}). Since $\lambda_n (n c^2_n)^{-1} \rightarrow \infty$ as $n \rightarrow \infty$, $\Gamma$ a compact set, thus  the second term of the right-hand side of  (\ref{ep}) converges to 0,  and then it is smaller than $\epsilon_3/2$ for $n$ large enough. Thus, the Lemma is proved by  taking $\epsilon_1=\epsilon_3/2$.\\
\hh Let us note that the adaptive weight $\hat \eo^t_{(0;n)}$  is calculated on all observations,  without imposing the constraint that $\| \ef -\efo\|_2 \geq c_n$.
\hspace*{\fill}$\blacksquare$ \\

  \noindent {\bf Proof of Lemma \ref{Lemma 9Bai}}
 By similar calculations to the proof of Proposition \ref{proposition 2.1.SCAD}, we obtain that
 \[
 \eE[R^{(\tau)}_i(b^0_r+\mu,\ef^0_r;b^0_r,\ef^0_r)] = \int^0_{- |\mu|} [x+\mu ] dF(x+b^0_r) 
 \]
is a positive function for any $\mu \in {\cal B}$ and a increasing function in $|\mu|$, with a single zero for $\mu=0$. By  Lemma \ref{Lemma 8Bai}, there exists at least $(l^0_r-l_r)\epsilon_0$ observations, for some $\epsilon_0 > 0$, and some $\delta> 0$,  such that $|b^0_{r+1}-b^0_r |+|\XX^t_i(\ef^0_{r+1})- \ef^0_r )| > \delta$. Then 
 \[
 \eE \cro{ \sum^{l^0_r}_{i=l_r+1}R^{(\tau)}_i(b^0_{r+1},\ef^0_{r+1};b^0_r,\ef^0_r) } \geq (l^0_r- l_r) \epsilon_0  \int^0_{- \delta} [x+\delta ] dF(x+b^0_r).
 \]
 By Lemma \ref{Lemma 4Bai}, for $\efo=\ef^0_r$, $\ef=\ef^0_{r+1}$ and $c_n=\max(|b^0_{r+1}-b^0_r |, \|\ef^0_{r+1})- \ef^0_r \|_2 )$, we have that for all $\epsilon >0$  such that 
 \[\eP [|\sum^{l^0_r}_{i=l_r+1}R^{(\tau)}_i(b^0_{r+1},\ef^0_{r+1};b^0_r,\ef^0_r) -\eE[\sum^{l^0_r}_{i=l_r+1}R^{(\tau)}_i(b^0_{r+1},\ef^0_{r+1};b^0_r,\ef^0_r) ] | > \epsilon(l^0_r-l_r) ] \]
 $ \leq  \exp(-C \epsilon(l^0_r-l_r)).$
  We take $\epsilon=2^{-1}\epsilon_0 \int^0_{- \delta} [x+\delta ] dF(x+b^0_r)$ and $\eta=2^{-1}\epsilon_0 \int^0_{- \delta} [x+\delta ] dF(x+b^0_r)$. Then  the  Lemma follows.
  \hspace*{\fill}$\blacksquare$ \\  
  
  \bibliographystyle{plain}
%\textbf{References}

\end{document}